\documentclass[11pt,a4paper,reqno]{amsart}
\usepackage{amssymb}
\usepackage{graphicx}
\usepackage{pdfsync}
\usepackage{xspace}
\usepackage[english]{babel}
\usepackage{amsfonts}
\usepackage{latexsym}
\usepackage{graphics}
\usepackage{hyperref}
\usepackage{pifont,bbding}
\usepackage{color}
\usepackage{bigints}
\usepackage{bm}
\usepackage{enumerate,amsthm,amsmath,comment,psfrag}
 \usepackage{mathrsfs} 
 \usepackage{nicefrac}

\newtheorem{lemma}{Lemma}[section]
\newtheorem{proposition}{Proposition}[section]
\newtheorem{theorem}{Theorem}[section]

\numberwithin{equation}{section}
\newcommand{\tnorm}[1]{\left\vert\!\left\vert\!\left\vert #1 
\right\vert\!\right\vert\!\right\vert}

\newcommand{\Rea}{\mathrm{Re}}
\newcommand{\Imag}{\mathrm{Im}}

\makeatletter
\newcommand*\bigcdot{\mathpalette\bigcdot@{.5}}
\newcommand*\bigcdot@[2]{\mathbin{\vcenter{\hbox{\scalebox{#2}{$\m@th#1\bullet$}}}}}
\makeatother

\title[A numerical scheme for the nonlinear von Neumann equation]
{
A second-order-in-time scheme for the von Neumann equation with singular self-interaction and  simulation of the onset of  instability}
\date{\today}

\author{Agissilaos  Athanassoulis}
\address[Agissilaos  Athanassoulis]{Department of Mathematics, University of Dundee, Dundee DD1 4HN, Scotland, UK} \email{a.athanassoulis@dundee.ac.uk}

\author{Fotini Karakatsani}
\address[Fotini Karakatsani]{Department of Mathematics, University of Ioannina, Ioannina 45110, Greece} \email{fkarakatsani@uoi.gr}

\author{Irene Kyza}
\address[Irene Kyza]{School of Mathematics and Statistics, University of St Andrews, St Andrews KY16 9AJ, Scotland, UK} \email{ika1@st-andrews.co.uk}

\begin{document}

\keywords{relaxation scheme, structure preserving, finite differences, Alber equation, modulation instability, nonlinear Landau damping.
}
\begin{abstract}
 The von Neumann equation with delta self-interaction kernel serves as a statistical model for nonlinear waves, and it exhibits a bifurcation between stable and unstable regimes. In oceanography it is known as the Alber equation, and its bifurcation is  important  for understanding rogue waves, a key problem in marine safety. Despite its significance, only one first-order-in-time numerical method exists in the literature.  In this paper, we propose a structure-preserving, linearly implicit, second-order-in-time scheme for its numerical solution. We employ fourth-order finite differences for the spatial discretization.   As an illustrative example, we explore the onset of modulation instability. We verify that the linear stability analysis accurately predicts the initial growth phase, but   fails to forecast the maximum amplitude, the formation of a coherent structure in the nonlinear regime, or the relevant  timescales. Monte Carlo simulations with Gaussian background spectra reveal that the maximum amplitude  depends mainly on the homogeneous background rather than the initial inhomogeneity. For weak instabilities, the inhomogeneity grows substantially from its initial condition, but  remains small compared to the background. On the other hand, strong instability leads to  recurrent hotspots of increased variance.  This provides a  possible explanation of how modulation instability makes rogue waves more likely in unidirectional sea states. 
\end{abstract}

\maketitle

\section{Introduction}

In  ocean waves,
the complex-valued envelope of a unidirectional sea state can be modeled with the focusing cubic nonlinear Schr\"odinger equation (NLS) \cite{mei2005theory},
\[
i \partial_t A + p \Delta A + \frac{q}2|A|^2 A = 0.
\]
The NLS is widely used as a starting point for the investigation of extreme events and rogue waves \cite{sapsis2021statistics,dematteis2018rogue}. 
Since power spectra are among the most widely available measurements of ocean wave fields, it is desirable to have a theory of NLS dynamics at the level of second moments \cite{dysthe2008oceanic,Ribal2013AlberEquation}.
In that context, the two-space autocorrelation is defined as
\[
R(x,y,t) = \mathrm{E} \big[ A(x,t) \overline{A(y,t)} \big]
\]
and, closing the BBGKY-type hierarchy at second order under a Gaussian approximation,
it is seen to satisfy
\[
i \partial_t R + p (\Delta_x - \Delta_y) R + q\Bigl(R(x,x,t)-R(y,y,t)\Bigr) \, R=0
\]
to leading order \cite{Alber1978}.  Empirically, it is well known that the second moment is typically close to homogeneous in appropriate temporal and spatial scales \cite{ochi1998ocean}. This insight can be made precise by the decomposition of the two-space autocorrelation as
\[
R(x,y,t) = \Gamma(x-y) + u(x,y,t),
\]
that is, into a homogeneous background $\Gamma$ and an inhomogeneous part $u,$ which is typically  much smaller than  $\Gamma.$ The inhomogeneity $u$  then satisfies the equation
\begin{equation}\label{eq:inhomalb1}
\begin{array}{c}
		i \partial_t u + p (\Delta_x - \Delta_y) u + q\Bigl(V(x,t)-V(y,t)\Bigr) \, \Bigl( \Gamma(x-y)  + u(x,y,t) \Bigr)=0,  \\[8pt]
		V(x,t) = u(x,x,t), \qquad u(x,y,0) = u_0(x,y),
\end{array}
\end{equation}
which is the focus of this paper. Equation \eqref{eq:inhomalb1} is called the Alber equation in oceanography. It can also be called a von Neumann equation with singular self-interaction to highlight its relationship with quantum statistical mechanics, discussed in Section \ref{sec:contextwig}.

Note that $p$ and $q$ are real-valued parameters, $\Gamma$ is a given autocorrelation function, and $u_0$ is Hermitian, $u_0(x,y) = \overline{u_0(y,x)}.$ By $\Delta_x$ we denote the Laplacian in the $x$ variable, while by $\Delta_y$ the Laplacian in the $y$ variable. The autocorrelation $\Gamma$ describes the known homogeneous background of the problem, and the unknown function $u$ represents the inhomogeneity. By the symmetry of the problem, $V(x,t)$ is always real valued.

\smallskip

  A question of particular interest is under what conditions nonlinear localized features can grow out of the quasi-homogeneous sea state — that is, under what conditions $u$ can grow to become comparable with $\Gamma.$  This question is closely related to the stability of the solution $u=0$ of \eqref{eq:inhomalb1}. 

Note that equation \eqref{eq:inhomalb1} contains the given background autocorrelation $\Gamma$ as part of the data of the  problem. Autocorrelations $\Gamma$ (or equivalently their Fourier transforms, i.e., power spectra) are the most widely used and widely available measured data about real-world sea states. They are also empirically known to be crucial for the qualitative behavior of the sea state. Indeed, equation \eqref{eq:inhomalb1}   exhibits a bifurcation between a stable regime (akin to nonlinear Landau damping) and an unstable regime (generalized modulation instability), depending on the homogeneous background $\Gamma$ \cite{Alber1978,Gramstad2017,Onorato2003LandauDamping,Ribal2013AlberEquation,Shukla2007ModulationalIncoherent}. The  implications of this bifurcation for ocean waves, and for extreme events in particular, have been explored by several authors  \cite{Athanassoulis2017RogueWaves,dysthe2008oceanic,onorato2009statistical,stiassnie2008recurrent}.

\smallskip 

Beyond its relevance for ocean waves, the well-posedness theory for equation \eqref{eq:inhomalb1} is  still developing. 
  Recent theoretical works on the well-posedness of equation \eqref{eq:inhomalb1} include \cite{Athanassoulis2018,Athanassoulis2023b,chen2017global,grekhneva2020dynamics,sakbaev2022blow}. It is important to note that for the focusing nonlinearities relevant to ocean waves (i.e. $pq>0$) there are no global-in-time existence results; only local-in-time solutions with a lower bound on the possible blow-up time \cite{Athanassoulis2018,chaub2025semiclassical,han2025semiclassical}. This is related to the loss of derivatives when controlling norms of the trace $V(x)=u(x,x)$ from norms of the two-dimensional function $u(x,y),$ and it is not clear whether  genuine blow-up typically occurs, or if the analysis is still not refined enough to preclude it under reasonable conditions. On the other hand, for  defocusing problems $(pq<0),$  there are some global-in-time existence results in the literature,  under additional conditions  \cite{chen2017global,hadama2025global}. 

\smallskip 

The main contributions of this paper are:
\begin{enumerate}
	\item The proposition, implementation and validation of scheme \eqref{eq:rCN1}-\eqref{eq:7bnvcx}. The need for a second-order initialization is documented and discussed in Section \ref{sec:2}.
	\item The preservation of invariants and discrete balance law of the scheme, discussed in detail in Section \ref{subsec:inv} for the time-discrete scheme and elaborated in Section \ref{sec:inarinants} for the fully discrete problem. 
	\item The consistency  and global boundedness of the scheme, proved in Sections \ref{sec:consistency} and \ref{sec:boundedness}.
	\item The numerical investigation of the fully nonlinear Landau damping and modulation instability in Section \ref{sec:LD}. The formation of a meta-stable pattern is found in the unstable regime. 
	\item Finally, in a Monte Carlo investigation of modulationally unstable sea states described in Section \ref{sec:4.3}, it is  seen that the behavior of the problem depends mainly on the background $\Gamma,$ and not on the exact initial condition.
Due to a quadratic dependence (cf. Figure \ref{fig:AFs2}), for barely unstable sea states the net physical effect is negligible, despite the presence of linear instability. For stronger instabilities the physical effect is substantial, and its implications for rogue waves   are discussed.
\end{enumerate}

\smallskip 

Before we go into more detail on the scheme, it is helpful to briefly go over existing numerical methods for the broader class of von Neumann equations, and for the Alber equation in particular.

\subsection{Broader context: von Neumann and Wigner equations}\label{sec:contextwig}

Equation \eqref{eq:inhomalb1} is in fact a special case of a broader class of von Neumann equations that appear in quantum statistical mechanics \cite{neumann1927wahrscheinlichkeitstheoretischer}. The von Neumann equation with background $\Gamma$ and  potential $V$ is
\begin{equation}\label{eq:vNstate}
	i \partial_t \rho + p (\Delta_x - \Delta_y) \rho + q(V(x,t)-V(y,t))(\Gamma(x-y) + \rho)=0, \qquad \rho(x,y,0)=\rho_0(x,y),
\end{equation}
and the potential $V$ may be a given function, or often 
it depends on the unknown function $\rho$ via
\begin{equation}\label{eq:nonlineariy}
	V(x,t) =  K(x) \ast \rho(x,x,t)
\end{equation}
for some self-interaction kernel $K.$ In the case of electrostatic self-interaction of a particle, $K$ would be  the Poisson kernel. If a smooth kernel $K$ is used, then equation \eqref{eq:vNstate} is often called a Hartree equation. In that context, the Alber equation amounts to taking a delta kernel, $K(x)=\delta(x)$ (hence, singular self-interaction). While an extensive literature exists  for the Hartree and Poisson cases dating back many decades,  the analysis for $K(x)=\delta(x)$ is much more recent and still developing \cite{Athanassoulis2018,chaub2025semiclassical,chen2017global,hadama2025global,han2025semiclassical}. Before we go into  the Alber equation in detail, it is worth briefly going over the context of numerical schemes for von Neumann-type equations at large.

In \cite{wigner1932quantum} it was pointed out that a particular transformation of the von Neumann equation \eqref{eq:vNstate} leads to a quantum kinetic-type equation (the Wigner equation), which in the semiclassical limit (an appropriate asymptotic scaling) formally converges to the kinetic equations of classical statistical mechanics. This provides a well known demonstration of the consistency between quantum and classical mechanics, which is still being refined, e.g. \cite{chong2024many}. The Wigner equation with its intuitive resemblance to kinetic equations has  been heavily used in several problems, including the modeling of resonant tunneling diodes \cite{Ringhofer1989AbsorbingBC,Arnold1994AbsorbingBC}.

A key point here is that there are many formally equivalent formulations of the same dynamics. In the von Neumann equation \eqref{eq:vNstate},  the independent variables $x,y$  both denote spatial variables. By applying the transform
\begin{equation}
W(x,k,t) \;=\; \mathcal{F}_{y \to k}\!\Bigl[\,
\rho\bigl(x+\tfrac{y}{2},\,x-\tfrac{y}{2},t\bigr)\Bigr]
\end{equation}
we pass to the Wigner equation on phase-space. Here $x$ is space and $k$ is wavenumber. The resulting equation reads
\begin{equation}\label{eq:wignereq1}
	\partial_t W+2 \pi p k \cdot \partial_x W - qT_V W=0 
\end{equation}
where the potential term $T_V W$ is \cite{lions1993mesures,athanassoulis2011strong}
\begin{align} \label{eq:TV}
T_V W =& i \int\limits_{\lambda, y \in \mathbb{R}^n} e^{-2 \pi i \lambda y}\left[V\left(x-\frac{y}{2}, t\right)-V\left(x+\frac{y}{2}, t\right)\right] dy \,W(x, k-\lambda, t) d \lambda  \\
=&  2 \operatorname{Re}\left[i \int\limits_{S\in\mathbb{R}^n} e^{2 \pi i S x} \widehat{V}(S) W^{\varepsilon}\left(x, k-\frac{ S}{2}\right) d S\right]  \\
=& \mathcal{F}_{X,K \to x,k}^{-1} \Bigl[ 
2 \int\limits_{S\in\mathbb{R}^n} \widehat{V}(S) \widehat{W}(X-S, K) \sin (\pi  S K) d S
\Bigr]
\end{align}

It is already evident that different formulations  may be advantageous for different purposes.  When designing a numerical method, the first decision concerns which formulation to discretize. While the von Neumann equation \eqref{eq:inhomalb1} and Wigner equation \eqref{eq:wignereq1} are formally equivalent at the level of the continuous problem, discretizing each one presents different trade-offs and challenges.

\medskip 


There has been substantial interest in the simulation of the von Neumann/Wigner dynamics  in general. Several works focus on Poisson nonlinearities and discretize  equation \eqref{eq:wignereq1}, typically using the potential term of \eqref{eq:TV} or its formal series expansion. This formulation allows  modelling of  real-world devices and impose boundary conditions at their endpoints modeling their actual connections, incoming current etc. However these same boundary conditions also cause difficulties, discussed below. In \cite{Jiang2023WignerPoissonRTD} a second-order-in-time method is introduced, with a sinc-Galerkin discretization for the potential term in the $k$ variable and finite differences in $x$. In \cite{Chen2022SpectralWignerPoisson,Chen2019UnboundedWigner,Chen2019Wigner4D}, fourth-order-in-time methods are proposed for the Wigner equation using the expression \eqref{eq:TV}. In \cite{Jiang2021HybridWigner} a second-order-in-time implicit-explicit method is introduced.
In \cite{Muscato2019StochasticWigner,Muscato2016StochasticWigner} stochastic algorithms are introduced. In \cite{Dorda2015WENO} a second-order-in-time method is proposed using a finite-volumes-like idea for the space discretization. Some classical works include a particle method \cite{arnold1992numerical}, and a first-order-in-time spectral collocation method \cite{ringhofer1992spectral}.  All of the aforementioned works discretize equation \eqref{eq:wignereq1}.

All of this activity highlights the fact that there is not an established family of  numerical methods for the Wigner-Poisson equation that stand out as the go-to option. This is due to some fundamental difficulties.
In one form or another, all of the aforementioned works truncate equation \eqref{eq:wignereq1} on a finite computational domain, and then proceed to discretize it. The main issue  is that the potential term $T_V W$ is global in  $(x,k)\in\mathbb{R}^{2d}.$ There is no well-posedness or propagation of regularity for \eqref{eq:wignereq1} on a finite computational domain, and that is why no systematic  analysis is possible. 
This difficulty has long been  known, and recently it has been  investigated more systematically \cite{rosati2013wigner}. One approach to bypass this issue is, for instance, to discretize in time only and analyze the resulting scheme at that level \cite{Arnold1996WignerPoisson}. 
More recently, another idea has been to avoid a hard truncation to a finite computational domain, and instead project on a Hermite basis of $L^2(\mathbb{R}^{2n})$ to implement the spatial operators; this was carried out for a linear problem in \cite{Filbet2025SemiClassicalVN}. 

Another way to work around the spatial truncation of a global term, is to focus on the two-space variables formulation \eqref{eq:vNstate} instead of the space-wavenumber formulation  \eqref{eq:wignereq1}. Applying e.g. periodic or Dirichlet boundary conditions on a spatially-truncated version of \eqref{eq:vNstate}  does not destroy well-posedness, even for strongly nonlinear problems. While the equation remains non-local (coupling the point $(x,y)$ to the points $(x,x)$ and $(y,y)$), there are no global terms invoking explicitly the  whole of space. 
Indeed, several authors have tried to simulate these dynamics by localizing and discretizing the von Neumann equation \eqref{eq:vNstate}, and subsequently taking Fourier transforms afterwards to pass to Wigner functions  if necessary.
In \cite{Schulz2021ExponentialIntegrators} the fundamental limitations of truncating and discretizing the Wigner equation \eqref{eq:wignereq1} are recognized, and a numerical method  involving an exponential integrator in time is proposed at the level of the von Neumann equation for a Hartree kernel. In \cite{Tian2012LvNFiniteTemp} a von Neumann-type equation with the Kohn-Sham Hamiltonian and dissipation is again discretized, while in \cite{Ribal2013AlberEquation,stiassnie2008recurrent} it is the Alber equation on rotated space variables (more details below).


\smallskip 

To summarize, while there are many numerical schemes for nonlinear von Neumann / Wigner equations,  almost all are concerned with Poisson or Hartree interaction kernels. There exists in the literature one numerical method explicitly designed and studied for the singular case of the Alber equation (i.e. $V(x,t):=u(x,x,t)$), and this is in the references \cite{Ribal2013AlberEquation,stiassnie2008recurrent}, further discussed below in  Section \ref{subsec:numeric}.

\subsection{Numerical methods for the von Neumann equation with singular interaction  (Alber equation): state of the art}\label{subsec:numeric}

As we saw in Section \ref{sec:contextwig}, the vast majority of existing numerical work for nonlinear von Neumann equations is for more regular interaction kernels. Here, we are focusing on the state of the art for the particular, singular nonlinearity that appears in the Alber equation.

The references \cite{Ribal2013AlberEquation,stiassnie2008recurrent}   treat equation  \eqref{eq:inhomalb1} after a change of variables, to ``mean position'' $X=\frac{x+y}2,$ and ``lag''   $r=x-y$ (this change of variables first appears in \cite{Alber1978}, and it is also employed in the calculations of Appendix \ref{app:A}). Both papers \cite{Ribal2013AlberEquation,stiassnie2008recurrent} use an explicit Euler  time-stepping scheme together with second-order finite differences in space. Moreover, a non-standard boundary condition  in the variable $r$ is used.

The works \cite{Ribal2013AlberEquation,stiassnie2008recurrent} are fundamental to the field, and they highlight the challenges this non-standard equation poses. The  boundary conditions in $r,$ introduced in \cite{stiassnie2008recurrent}, attempt to match the far-field behavior of the equation with a known background behavior derived from a box spectrum.  In \cite{Ribal2013AlberEquation} the boundary condition is modified to match a JONSWAP spectrum. 
The underlying rationale for these boundary conditions  is that the inhomogeneity should vanish at infinity; they represent a natural attempt to enforce  the ``zero at large $r$'' constraint.
 As we see in the numerical results, in the case of modulationally unstable problems, the inhomogeneity extends its support, and inevitably reaches the boundary of the computational domain. 

Our approach is to solve directly \eqref{eq:inhomalb1} with periodic boundary conditions in the $x,y$ variables. This preserves by construction the original physical context \cite{ochi1998ocean}, and is known to produce equivalent results to whole space for large enough computational domain \cite{Athanassoulis2023b}.  This also treats in a unified way any background spectrum (i.e. there is no need to specially fine tune the boundary conditions to each spectrum as in \cite{Ribal2013AlberEquation}).  In any case, the high order of convergence of our method allows us to use sufficiently large computational domains so that the modulation instability is fully developed before any interaction with the boundary (cf., e.g.,  Figure \ref{fig:3}).

Crucially, our scheme  satisfies a discrete version of the  balance law controlling the growth of the solution  (see Propositions \ref{lm:msbal21} and  \ref{prop:discrconbll}), as well as preserving  two invariants at discrete level (see Section \ref{sec:inarinants}). This is particularly important in unstable cases, where the solution $u$ grows by a factor of tens or hundreds.

\subsection{Time semi-discrete scheme} The scheme we propose for the simulation of equation \eqref{eq:inhomalb1} can be described as a relaxation--Crank--Nicolson in time, and is inspired by \cite{Athanassoulis2023a,besse2004relaxation,besse2021energy}.
We start from equation \eqref{eq:inhomalb1} and write it in system form using the auxiliary variable $\phi$
\begin{equation}\label{eq:aug1}
\begin{array}{c}
\displaystyle  \phi(x,y,t) = u(x,x,t) - u(y,y,t), \\[4pt]
\displaystyle  i \partial_t u(x,y,t) + p (\Delta_x - \Delta_y) u(x,y,t) + q \phi(x,y,t) \,  \Bigl(\Gamma(x-y) + u(x,y,t) \Bigr)  = 0 \\
\end{array}
\end{equation}
on $[-\frac{L}2,\frac{L}2]$ equipped with periodic boundary conditions, 
\[
 (x,y) \in [-\frac{L}2,\frac{L}2]^2, \quad \partial_y^j u(x,-\frac{L}2) = \partial_y^j u(x,\frac{L}2), \,\, \partial_x^j u(-\frac{L}2,y) = \partial_x^j u(\frac{L}2,y), \quad j\in\{0,1\}.
\]
Note that $\Gamma$ must also be an $L$-periodic function. In practice this means that the size $L$ of the domain must be large enough so that  the effective support of $\Gamma$ is contained in $[-\frac{L}2,\frac{L}2].$

Now equation \eqref{eq:aug1} is discretized in time to
\begin{equation}\label{eq:rCN1}
\begin{array}{c}
\displaystyle  \frac{1}2 \left(  \Phi^{n+\frac{1}2}(x,y) + \Phi^{n-\frac{1}2}(x,y) \right) = U^n(x,x) - U^n(y,y), \\[9pt]
\displaystyle  i \frac{ U^{n+1} - U^n }{\tau}+ p (\Delta_x - \Delta_y) U^{n+\frac{1}2} +  q \Phi^{n+\frac{1}2} \, \Bigl(\Gamma(x-y) + U^{n+\frac{1}2} \Bigr) = 0, \\
\end{array}
\end{equation}
where
\begin{equation}\label{eq:7bnvcx}
	U^{n+\frac{1}2}=U^{n+\frac{1}2}(x,y) := \frac{U^{n+1}(x,y)+U^{n}(x,y)}2.
\end{equation}

This inherits the main symmetries of the original problem; in particular $\Phi^{n+\frac{1}2}$ is always real-valued.

Note that the first equation in \eqref{eq:rCN1} is obtained by a linear extrapolation from the points $t^n$ and $t^{n-\frac{1}2}.$ 
Moreover, the second equation in \eqref{eq:rCN1} is, essentially, a Crank--Nicolson time discretization of the second equation in \eqref{eq:aug1}. 
As expected, this leads to a second-order-in-time scheme. 
The auxiliary variable $\Phi^{n+\frac{1}2}$ lives on a staggered time grid $t^{n+\frac{1}2}:=(t^{n+1}+t^n)/2.$ Its introduction removes the need for solving a nonlinear equation in each timestep, while keeping strong stability features -- analogous to those of the fully implicit Crank--Nicolson. The price for that is requiring  the additional initialization of $\Phi^{-\frac{1}2}$ for the first timestep. The initialization is discussed in more detail in Section \ref{sec:initimpl}.

\medskip

The remainder of the paper is structured as follows:
Sections \ref{sec:2m} and \ref{sec:2} provide a self-contained description of the theory for the periodized problem \eqref{eq:aug1}, as well as for the  numerical method, its features and validation. In Section \ref{sec:LD} we use the numerical method to investigate key questions about equation \eqref{eq:inhomalb1} (or, more precisely, its periodized version \eqref{eq:aug1}), such as the onset of instability, the maxima of the fully developed instability, and the dispersion of inhomogeneities in the stable case. By necessity, Section \ref{sec:LD} relies heavily on the state of the art for equation \eqref{eq:inhomalb1}. 
%

\section{Theoretical results} \label{sec:2m}

\subsection{Well-posedness of the continuous problem}

It must be emphasized that the Alber equation contains a singular nonlinearity, and most of the standard tools used in the analysis of von Neumann / Wigner equations no longer apply. In particular, there are no global-in-time existence results for the Alber equation with focusing nonlinearity $(p q>0)$ in the literature. 
The state of the art  provides local-in-time existence and propagation of regularity in space and time (up to a possible blow-up time) \cite{Athanassoulis2018}. The main difficulties are the trace nonlinearity, and the absence of any a priori bounds (energy etc).

Problem \eqref{eq:aug1} is the Alber equation localized on the $L\times L$ torus $\mathbb{T}_L$. Regularity follows analogously to the full-space problem $(x,y)\in \mathbb{R}^{2d}$,  by adapting the Fourier-transform-based arguments to the Fourier series
\begin{equation}\label{eq:fcoeffs}
	u(x,y,t) = \sum\limits_{k,l \in \mathbb{Z}} \hat{u}_{k,l}(t) e^{2\pi i \frac{kx+ly}L}, \qquad 
	\Gamma(y) = \sum\limits_{n\in\mathbb{Z}} P_n e^{2\pi i \frac{ny}L}.
\end{equation}

\begin{theorem}\label{thrm:torusreg}
Consider problem \eqref{eq:aug1} with Hermitian initial data $u(0) = u_0$ satisfying 
\[
\tnorm{u_0} := \sum\limits_{k,l} |\hat{u}_{k,l}(0)| <\infty.
\] 
Then, there exists a time $T_*$ depending only on $q, \Gamma$ and the size of the initial condition $\tnorm{u_0}$ such that problem \eqref{eq:aug1} has a unique solution $u(t),$ with  $\tnorm{u(t)} <\infty$ for $t\in [0,T_*).$

The solution propagates regularity in the sense that
\[
\max\limits_{|a+b|\leqslant s} \tnorm{ \partial_x^a \partial_y^b u_0} <\infty \implies 
\max\limits_{|a+b|\leqslant s} \tnorm{ \partial_x^a \partial_y^b u(t)} <\infty \quad \forall t\in [0,T_*),
\]
and smoothness in time can be shown in the sense that
\[
\max\limits_{|a+b|\leqslant 2m} \tnorm{ \partial_x^a \partial_y^b u_0} <\infty \implies 
 u(t) \in C^m [0,T_*).
\]
\end{theorem}

The proof follows the logic of \cite{Athanassoulis2018}, and for completeness it is included in Appendix \ref{App:A}.

Observe that $\max\limits_{|a+b|\leqslant s} \tnorm{ \partial_x^a \partial_y^b u(t)} <\infty$ means that $u(t) \in H^s(\mathbb{T}_L).$

\subsection{Balance laws and invariants} \label{subsec:inv}

An $L^2$ norm balance law for the general continuous problem relates the rate of growth of the solution with the data of the problem. Our numerical scheme inherits this balance law at the discrete level.

On the continuous level, the balance law is as follows:
\begin{proposition}\label{lm:msbal21} Consider equation \eqref{eq:inhomalb1}. Assume that the initial condition $u_0$ is a smooth function, and that the given autocorrelation function $\Gamma(x)$ is of the form $\Gamma(x) = \mathcal{F}^{-1}_{k \to x}[P(k)],$ where $P(k)$ is a nonnegative smooth function of compact support.

Then the potential $V(x,t)=u(x,x,t)$ is real valued for all times.
Moreover, if we denote
$\mathcal{M}(t)  := \|u(t)\|_{L^2}^2,$
it follows that 
\begin{equation}\label{eq:msbal21}
\frac{d}{dt}\mathcal{M}(t) = 2 \Rea \Bigl[     iq \int\limits_{x,y} \Bigl(V(x,t)-V(y,t)\Bigr)\Gamma(x-y) \overline{u(x,y,t)} dxdy  \Bigr].	
\end{equation}
Equation \eqref{eq:msbal21} holds when $(x,y) \in \mathbb{R}^{2}$ and $u_0,$ $\Gamma$ exhibit rapid decay, or when $(x,y) \in [-\frac{L}2,\frac{L}2]^{2}$ and $u,$ $\Gamma$ satisfy periodic boundary conditions.
\end{proposition}

\noindent {\bf Proof:} Recall that the initial condition $u_0(x,y)$ for \eqref{eq:inhomalb1} is Hermitian, $u_0(x,y) = \overline{ u_0 (y,x)}.$ This Hermitianity  is preserved in time, 
and implies that $V(x,t)=u(x,x,t)$ is real-valued for all times. The real-valuedness of the potential is crucial in what follows.
Because of the well-posedness results discussed above,
 all the integrals appearing below exist.

With regard to the growth of $\mathcal{M}(t),$ we readily observe that
\[
\begin{aligned}
	\frac{d}{dt} \mathcal{M}(t) & =\frac{d}{dt} \langle u, u \rangle =  2 \Rea \Bigl[  \langle u_t, u \rangle  \Bigr] 	 \\
	& =2 \Rea \Bigl[  \Bigl\langle i p (\Delta_x - \Delta_y)u + iq \Bigl( V(x,t)-V(y,t) \Bigr) \, \Bigl( \Gamma(x-y)+ u \Bigr) \, , \,\, u \Bigr\rangle  \Bigr] \\
	& =2 \Rea \Bigl[  -ip \Bigl\langle  \nabla_x u,\nabla_x u \Bigr\rangle+ip \Bigl\langle  \nabla_y u,\nabla_y u \Bigr\rangle   +  \hfill \\
& \qquad\qquad	 +\Bigl\langle iq \Bigl( V(x,t)-V(y,t) \Bigr) \, \Bigl( \Gamma(x-y)+ u \Bigr) \, , \,\, u \Bigr\rangle  \Bigr] \\
& =2 \Rea \Bigl[     iq \int\limits_{x,y} \Bigl(V(x,t)-V(y,t)\Bigr)\Gamma(x-y) \overline{u(x,y,t)} dxdy  \Bigr] .
\end{aligned}
\]
The boundary conditions come into play only at the integration by parts of the Laplacian terms; the proof applies in the same way for full space $\mathbb{R}^{2d}$ with decay at infinity or if $u$ satisfies periodic boundary conditions. \qed

\medskip

Essentially the same computation works out at the discrete level for the numerical scheme \eqref{eq:rCN1}. We present it here in the time-discrete form for simplicity. The same carries over to the fully discrete level, where the discrete $L^2$ norm would appear in the place of the continuous $L^2$ norm.

\medskip


\begin{proposition}\label{prop:discrconbll}
Consider the time-discrete scheme \eqref{eq:rCN1}, and denote
\begin{equation}\label{eq:dcqw21}
\mathcal{M}^n := \| U^n\|_{L^2}^2.
\end{equation}
Then the following time-discrete balance law holds:
\begin{equation}\label{eq:7y6t5r4eghtm}
\frac{\mathcal{M}^{n+1} - \mathcal{M}^n}{\tau} = 2 \Rea \Bigl[ i  q
\Bigl\langle \Phi^{n+\frac{1}2} \Gamma(x-y), U^{n+\frac{1}2} \Bigr\rangle
\Bigr].
\end{equation}
\end{proposition}

\noindent {\bf Remark:} Equation \eqref{eq:7y6t5r4eghtm} is a time-discrete version of equation \eqref{eq:msbal21}. 

\noindent {\bf Proof:} 
Observe that
\[
\begin{aligned}
2 \Rea \left[\left\langle U^{n+1}-U^n,\, U^{n+\frac12} \right\rangle\right]
&= \left\langle U^{n+1}-U^n,\, U^{n+\frac12} \right\rangle
 + \left\langle U^{n+\frac12},\, U^{n+1}-U^n \right\rangle \\
&= \frac12 \left\langle U^{n+1}-U^n,\, U^{n+1}+U^n \right\rangle
 + \frac12 \left\langle U^{n+1}+U^n,\, U^{n+1}-U^n \right\rangle \\
&= \|U^{n+1}\|_{L^2}^2 - \|U^n\|_{L^2}^2 .
\end{aligned}
\]
Thus,
\[
\begin{aligned}
\mathcal{M}^{n+1}-\mathcal{M}^n
&=2\,\Rea\!\left[
\langle U^{n+1}-U^n,\,U^{n+\frac12}\rangle
\right] \\
&=2\,\Rea\!\left[
\left\langle
i\tau p(\Delta_x-\Delta_y)U^{n+\frac12}
+i\tau q\,\Phi^{n+\frac12}(\Gamma(x-y)+U^{n+\frac12}),
\,U^{n+\frac12}
\right\rangle
\right] \\
&=2\,\Rea\!\Bigl[
i\tau p\langle \Delta_xU^{n+\frac12},U^{n+\frac12}\rangle
-i\tau p\langle \Delta_yU^{n+\frac12},U^{n+\frac12}\rangle \\
&\qquad
+i\tau q\langle \Phi^{n+\frac12}\Gamma(x-y),U^{n+\frac12}\rangle
+i\tau q\langle \Phi^{n+\frac12}U^{n+\frac12},U^{n+\frac12}\rangle
\Bigr] \\
&=2\,\Rea\!\left[
i\tau q\,
\langle \Phi^{n+\frac12}\Gamma(x-y),U^{n+\frac12}\rangle
\right].
\end{aligned}
\]\qed 

It should be noted that  the continuous problem has several invariants, namely
\begin{equation} \label{eq:i0}
I_0[u(t)] := \int\limits_{x,y} |\Gamma(x-y) + u(x,y,t)|^2 dxdy,
\end{equation}
\begin{equation}\label{eq:i1}
I_1[u(t)] := \int\limits_x u(x,y,t)\Big|_{y=x} dx , 	
\end{equation}
\begin{equation}
I_2[u(t)] :=  \int\limits_x (\partial_x - \partial_y)u(x,y,t)\Big|_{y=x} dx , 
\end{equation}
\begin{equation}
I_3[u(t)] := \left( \frac{q}p \int\limits_x u^2(x,y,t)\Big|_{y=x} dx + \int\limits_x (\partial_x - \partial_y)^2u(x,y,t)\Big|_{y=x} dx \right),
\end{equation}
i.e., $\frac{d}{dt} I_j[u(t)]=0$ for $j\in \{0,1,2,3\}.$ Invariant $I_0$ coincides with the quantity denoted by $\mathcal{N}(t)$ in Proposition \ref{prop:bddns}, and a discrete version of it is exactly preserved at the fully discrete level (cf. Section \ref{sec:inarinants}).
Invariants $I_j,$ $j\in\{1,2,3\}$ first appeared in \cite{stiassnie2008recurrent} and are derived in our notation for completeness in Lemma \ref{lm:inva} in the Appendix.  They are based on the symmetries of $u,$ and stay constant even if $u$ grows by factors of hundreds, as in the modulationally unstable examples we will see later (i.e. they involve large cancellations in general). Invariant $I_1$ is also preserved at the fully discrete level, cf. Section \ref{sec:inarinants}.
We will monitor the numerical behavior of the discrete invariants $I_2,I_3$ as a diagnostic for the quality of the numerical solution. In general $I_2$ can be preserved reasonably well even when the solution grows several orders of magnitude. On the other hand, $I_3$ involves second derivatives, and can only be controlled when the solution does not grow dramatically.

\subsection{Consistency}\label{sec:consistency}
In this section, we compute the consistency  error for the time-discrete scheme \eqref{eq:rCN1}.  The order of the consistency error in time is the same for the fully discrete scheme below, and can be obtained using similar arguments. To this end, we define the residuals $r_1^n(x,y)$ and $r_2^n(x,y)$ as the amount by which the exact solution $(\phi,u)$ of the system \eqref{eq:aug1} fails to satisfy the system \eqref{eq:rCN1} of the numerical scheme. Thus,
\begin{equation}
\label{residuals}
\begin{aligned}
&r_1^n(x,y):=\frac12\left(\phi(x,y,t^{n+\frac12})+\phi(x,y,t^{n-\frac12})\right) - u(x,x,t^n)+u(y,y,t^n)\\
&r_2^n(x,y):= i\dfrac{u(x,y,t^{n+1})-u(x,y,t^n)}{\tau}+p(\Delta_x-\Delta_y)u(x,y,t^{n+\frac12})+\\&\hspace{1.75cm}+q\phi(x,y,t^{n+\frac12})\left(\Gamma(x-y)+u(x,y,t^{n+\frac12})\right).
\end{aligned}
\end{equation}
To compute the order of the truncation error, we assume sufficient regularity on the solution $(\phi,u)$ of  \eqref{eq:aug1} -- an assumption that is based on Theorem \ref{thrm:torusreg} for sufficiently large $s,m.$ We also avoid writing the dependence on $(x,y)$, unless it is not clear. Since in this setting the equations have not been multiplied by $\tau,$ the expected order of the scheme in $\tau$ is the same as the order of $r^n_1,$ $r^n_2.$

We start by estimating $r_1^n$. We first note that
\begin{equation*}
\begin{aligned}
\phi(x,y,t^{n+\frac12})+\phi(x,y,t^{n-\frac12}) &=\phi(t^n)+\dfrac{\tau}2\partial_t\phi(t^n)+\dfrac{\tau^2}8\partial_{tt}\phi(t^n)+\dfrac{\tau^3}{48}\partial_{ttt}\phi(t^n)+\mathcal{O}(\tau^4)\\
&\quad +\phi(t^n)-\dfrac{\tau}2\partial_t\phi(t^n)+\dfrac{\tau^2}8\partial_{tt}\phi(t^n)-\dfrac{\tau^3}{48}\partial_{ttt}\phi(t^n)+\mathcal{O}(\tau^4)\\
&=2\phi(t^n)+\dfrac{\tau^2}4\partial_{tt}\phi(t^n)+\mathcal{O}(\tau^4).
\end{aligned}
\end{equation*}
Hence, using the above relation and the first equation in \eqref{eq:aug1},  we obtain
\begin{equation}
\label{firstresidual}
r_1^n=\frac {\tau^2}{8} \partial_{tt}\phi(t^n)+\mathcal{O}(\tau^4)=\mathcal{O}(\tau^2).
\end{equation}

We next estimate $r_2^n$. Since we use the first order forward approximation for the first derivative, we get
\begin{equation}
\label{auxr21}
i\dfrac{u(t^{n+1})-u(t^n)}{\tau}=i\left(\partial_t u(t^n)+\dfrac{\tau}2\partial_{tt}u(t^n)+\mathcal{O}(\tau^2)\right).
\end{equation}
Furthermore,
\begin{equation}
\label{auxr22}
p(\Delta_x-\Delta_y)u(t^{n+\frac12})=p(\Delta_x-\Delta_y)\left(u(t^n)+\dfrac{\tau}2\partial_tu(t^n)\right)+\mathcal{O}(\tau^2).
\end{equation}
In addition,
\begin{equation}
\label{auxr23}
\begin{aligned}
q\phi(t^{n+\frac12})  \left(\Gamma(x-y)+u(t^{n+\frac12})\right)
&= q\left(\phi(t^n)+\frac{\tau}2\partial_t\phi(t^n)+\mathcal{O}(\tau^2)\right)\\
&\quad \times\left(\Gamma(x-y)+u(t^n)+\frac{\tau}2\partial_tu(t^n)+\mathcal{O}(\tau^2)\right)\\
&=q\phi(t^n)\left(\Gamma(x-y)+u(t^n)\right)+q\dfrac{\tau}2\phi(t^n)\partial_tu(t^n)\\ 
&\quad +q\dfrac{\tau}2\partial_t\phi(t^n)\left(\Gamma(x-y)+u(t^n)\right)+\mathcal{O}(\tau^2).
\end{aligned}
\end{equation}
Combining \eqref{auxr21}, \eqref{auxr22} and \eqref{auxr23}, and using the second equation of \eqref{eq:aug1} (as well as taking its time-derivative) we get that
\begin{equation*}
\begin{aligned}
r_2^n=&i\left(\partial_tu(t^n)+\dfrac{\tau}2\partial_{tt}u(t^n)\right) +p(\Delta_x-\Delta_y)\left(u(t^n)+\dfrac{\tau}2\partial_tu(t^n)\right)+q\phi(t^n)\left(\Gamma(x-y)+u(t^n)\right)\\
&\quad+q\dfrac{\tau}2\phi(t^n)\partial_tu(t^n)+q\dfrac{\tau}2\partial_t\phi(t^n)\left(\Gamma(x-y)+u(t^n)\right)+\mathcal{O}(\tau^2)\\
=& i\partial_t u(t^n)+p(\Delta_x-\Delta_y)u(t^n)+q\phi(t^n)\left(\Gamma(x-y)+u(t^n)\right)\\
&\quad +\frac{\tau}2\bigl(i\partial_{tt}u(t^n)+q\partial_t\left(\phi(t)\left(\Gamma(x-y)+u(t)\right)\right)|_{t=t^n}+p(\Delta_x-\Delta_y)\partial_tu(t^n)\bigr)+\mathcal{O}(\tau^2)\\
=&\mathcal{O}(\tau^2).
\end{aligned}
\end{equation*}
Hence, we have shown that
$$r_1^n=\mathcal{O}(\tau^2),\,\, n\ge 1,\quad  \text{ and }\quad  r_2^n=\mathcal{O}(\tau^2), \,\, n\ge 0.$$
For $n=0$, we have $\Phi^{-\frac12}$ in the first equation of the time-discrete scheme. Therefore the order of  $r_1^0$ depends on the initialization of $\Phi^{-\frac12}$. If we are only interested in the order of $u$,  a first order initialization is fine.   However,  a reduction of order in $\phi$ is observed if we use the naive first-order initialization. To address this,  an advanced, second-order initialization is proposed.  This issue is discussed in detail in Sections \ref{sec:initimpl}, \ref{sec:EOC}.

It should be noted that, for the standard cubic nonlinear Schr\"odinger equation (NLS), there is no loss of order for the auxiliary variable with first-order initialization   \cite{besse2004relaxation, zouraris2023error}.
However, loss of order in the auxiliary variable with first-order initialization has been reported for analogous relaxation schemes in other problems \cite{athanassoulis2024efficient,Athanassoulis2023a,zouraris2021error}. More recently, advanced second order initializations are being widely used for this kind of relaxation scheme, even for the NLS \cite{ besse2021energy}.

\subsection{Boundedness and a remark on stability}\label{sec:boundedness}

There exists an unconditional $L^2$ boundedness result for $u$ that is inherited by $U^n:$ 

\begin{proposition}\label{prop:bddns} Consider the problem \eqref{eq:aug1} and denote
\[
\mathcal{N}(t) := \int\limits_{x,y \in \mathbb{T}_L} |\Gamma(x-y) + u(x,y,t)|^2 dxdy.
\]
Then
$
\frac{d}{dt}
\mathcal{N}(t) = 0 
$	
and there exists a constant $M,$ depending only on $u_0,$ $\Gamma$ and $L$ so that
\[
\| u(t)\|_{L^2} \leqslant M.
\]
Moreover, let $U^n$ be defined as in \eqref{eq:rCN1}, and denote
\[
\mathcal{N}^n :=  \int\limits_{x,y \in \mathbb{T}_L} |\Gamma(x-y) + U^n(x,y)|^2 dxdy.
\]
Then
$
\mathcal{N}^{n+1} = \mathcal{N}^n
$
and there exists a constant $M',$ depending only on $u_0,$ $\Gamma,$ and $L$ so that
\[
\| U^n\|_{L^2} \leqslant M'
\]
for any choice of timestep $\tau.$
\end{proposition}

\noindent {\bf Proof:} For the continuous problem, as long as the solution exists, we have
\[
\begin{aligned}
	\frac{d}{dt} \langle \Gamma(x-y) + u(t),&\Gamma(x-y) + u(t) \rangle = 2\Rea \left[ \langle \partial_t u(t)  , \Gamma(x-y) + u(t)\rangle \right]\\[5pt]
	&= 2\Rea \left[ \langle ip (\Delta_x - \Delta_y) u(t)  + i q \phi(\Gamma(x-y) + u(t))  , \Gamma(x-y) + u(t)\rangle \right] \\[5pt]
	&= 2\Rea \left[ \langle ip (\Delta_x - \Delta_y) (\Gamma(x-y)+ u(t))  , \Gamma(x-y) + u(t)\rangle  \right. \hfill\\  
	&\qquad\qquad +\hfill\left.  \langle  i q \phi(\Gamma(x-y) + u(t))  , \Gamma(x-y) + u(t)\rangle \right]=0
\end{aligned}
\]
Both terms are purely imaginary, therefore the whole expression vanishes. In the second to last step we used the fact that $(\Delta_x - \Delta_y) \Gamma(x-y) =0.$ We  thus conclude that $\displaystyle \frac{d}{dt}\mathcal{N}(t) = 0,$ which implies
\[
 \| \Gamma(x-y) + u(t)\|_{L^2} = \|\Gamma(x-y) + u_0\|_{L^2}.
\]
Hence it follows that
\[
\|u(t)\|_{L^2} 
\leqslant  \|\Gamma(x-y) + u_0\|_{L^2} + \|\Gamma(x-y)\|_{L^2}.
\]

All the steps above are essentially reproduced in the time-discrete case:
\[
\begin{aligned}
\mathcal{N}^{n+1}-\mathcal{N}^n
&=2\Rea \left[
\langle U^{n+1}-U^n,\Gamma(x-y)+U^{n+\frac12}\rangle
\right] \\
&=2\Rea \left[
\left\langle i\tau p(\Delta_x-\Delta_y)U^{n+\frac12}
+i\tau q\Phi^{n+\frac12}(\Gamma(x-y)+U^{n+\frac12}),
\Gamma(x-y)+U^{n+\frac12}
\right\rangle
\right] \\
&=2\Rea \left[
\left\langle i\tau p(\Delta_x-\Delta_y)(\Gamma(x-y)+U^{n+\frac12}),
\Gamma(x-y)+U^{n+\frac12}
\right\rangle \right.\\
&\qquad\qquad\left.
+\left\langle i\tau q\Phi^{n+\frac12}(\Gamma(x-y)+U^{n+\frac12}),
\Gamma(x-y)+U^{n+\frac12}
\right\rangle
\right]
=0 .
\end{aligned}
\]
Thus, similarly to what we had above,
\[
\| \Gamma(x-y) + U^n\|_{L^2} \leqslant \| \Gamma(x-y) + U^0\|_{L^2} + \| \Gamma(x-y)\|_{L^2}.
\]
Note that the constants $M,M'$ and the values $\mathcal{N}(t),$ $\mathcal{N}^0$ agree up to the approximation error of the initial condition, $\|u_0 - U^0\|_{L^2}.$
\qed

\smallskip

There are classical theoretical results for the Crank--Nicolson method applied to  Schr\"odinger-type equations where boundedness (as in Proposition \ref{prop:bddns}) leads to stability.  One may ask whether a similar argument applies in the current setting. Here, we want to briefly discuss the difficulty with this.

Consider the system \eqref{eq:rCN1} with initial conditions $U^0, \Phi^{-\frac{1}2},$ and the perturbation $\widetilde{U}^n,$ $\widetilde{\Phi}^{n+\frac{1}2}$ resulting from the perturbed initial conditions $\widetilde{U}^0, \widetilde{\Phi}^{-\frac{1}2}.$ Denote the errors by
\begin{equation}
	h^n := U^n - \widetilde{U}^n, \qquad f^{n+\frac{1}2} := \Phi^{n+\frac{1}2} - \widetilde{\Phi}^{n+\frac{1}2}
\end{equation}
(note that $f^{n+\frac{1}2}$ is real valued since $\Phi^{n+\frac{1}2},$ $\widetilde{\Phi}^{n+\frac{1}2})$ are real-valued).
Then the error equations are
\begin{equation}\label{eq:rCN1errpreq}
\begin{array}{c}
\displaystyle  f^{n+\frac{1}2}(x,y) = 2\big(h^n(x,x) - h^n(y,y) \big) - f^{n-\frac{1}2}(x,y), \\[9pt]
\displaystyle  
\frac{ h^{n+1} - h^n}\tau =  i   p (\Delta_x - \Delta_y) h^{n+\frac{1}2} + i  q f^{n+\frac{1}2} \, \Bigl(\Gamma(x-y) + U^{n+\frac{1}2} \Bigr) 
 +   i  q \widetilde{\Phi}^{n+\frac{1}2} \, h^{n+\frac{1}2} 
\end{array}
\end{equation}
where of course $h^{n+\frac{1}2}:=(h^{n+1}+h^n)/2.$ Since
\[
\| h^{n+1}\|^2_{L^2} - \|h^n\|^2_{L^2}= 2\Rea \left[\langle h^{n+1} - h^n, h^{n+\frac{1}2}\rangle \right],
\]
by using \eqref{eq:rCN1errpreq} to express $h^{n+1} - h^n$ we obtain
\[
\| h^{n+1}\|^2_{L^2} - \|h^n\|^2_{L^2}=2\tau\Rea \left[\langle    i  q f^{n+\frac{1}2} \, \Bigl(\Gamma(x-y) + {U}^{n+\frac{1}2} \Bigr)  , h^{n+\frac{1}2}\rangle \right]
\]
where the Laplacian and the symmetric potential term are cancelled out by the antisymmetry. This takes full advantage of the $L^2$-unitarity of the continuous problem's propagator, just like in the classical argument.
One obvious remaining difficulty is that the background $\Gamma$ can inject energy in the system, enlarging the solution (as in modulation instability) and also the error with it. This could likely be treated in the stable case, but the nonlinear term still could not. Even if we set $\Gamma=0$ for simplicity, we would get an estimate of the form
\begin{equation}\label{eq:nobrmgfr}
\| h^{n+1}\|^2_{L^2} - \|h^n\|^2_{L^2}  \leqslant 2q\tau\Rea \left[\langle   i   f^{n+\frac{1}2}  {U}^{n+\frac{1}2} , h^{n+\frac{1}2}\rangle \right].
\end{equation}
The natural norms here are $\|f^{n+\frac{1}2}\|_{L^\infty},$ $\|\widetilde{U}^{n+\frac{1}2} \|_{L^2}$ and  $\|h^{n+\frac{1}2}\|_{L^2}.$ The bottleneck is $\|f^{n+\frac{1}2}\|_{L^\infty};$ indeed, there is no clear way how to propagate the $L^\infty$ norm with
\begin{equation}\label{eq:trmnbvcffd}
 f^{n+\frac{1}2}(x,y) = 2\big(h^n(x,x) - h^n(y,y) \big) - f^{n-\frac{1}2}(x,y),
\end{equation}
when we use the $L^2$ norm for $h^n.$
The trace is not a bounded operator in $L^2,$ and there are no a priori bounds available at this point. Thus,  at the moment, we see no way to use equation \eqref{eq:nobrmgfr} to setup a discrete Gronwall inequality for $\|h^n\|_{L^2}$ as in the classical argument.

\section{Numerical validation}\label{sec:2}

\subsection{Fully discrete scheme}

To implement the numerical scheme, finite differences are used on a uniform $N\times N$ mesh  on $[-\frac{L}2,\frac{L}2]^2$ with spatial mesh size $h,$ $N h = L.$ Note that this  choice of mesh facilitates the implementation of $U(x,x,t)$ (which  is essentially a trace operator) and preserves the symmetry in $x$ and $y.$  Let  ${U}_{i,j}^{n}$  and ${\Phi}_{i,j}^{n+\frac{1}2}$ denote  the arrays of values of $U^{n}$  and $\Phi^{n+\frac{1}2}$  at the mesh points $(x_i,y_j)=(-\frac L2 + ih,-\frac L2 + jh), i, j=0, \ldots, N-1,$ respectively.

Fourth-order central differences with periodic boundary conditions are used for the spatial discretization. The Kronecker product is used to create a sparse $N^2 \times N^2$ matrix $\mathfrak{D}_H$ implementing the hyperbolic Laplacian $\Delta_x-\Delta_y$ on a vectorized form of $U^{n+\frac{1}2}_{i,j}$. Moreover, a $N^2\times N^2$ diagonal matrix $\mathfrak{F}^{n+\frac{1}2} = \texttt{diag}(\texttt{vec}({\Phi}_{i,j}^{n+\frac{1}2}))$ is created to implement the pointwise product of $\Phi^{n+\frac{1}2}$ with $U^{n+\frac{1}2}.$  The timestep  of equation \eqref{eq:rCN1} is executed by first updating $\Phi^{n+\frac{1}2}_{i,j},$ then solving for $U^{n+\frac{1}2}_{i,j},$ and finally extracting $U^{n+{1}}_{i,j}$ from it, as is standard for Crank-Nicolson-type schemes:
\begin{align}
&	\Phi^{n+\frac{1}2}_{i,j} = 2\Bigl( U^n_{i,i} - U^n_{j,j} \Bigr) - 	\Phi^{n-\frac{1}2}_{i,j}, \label{eq:23} \\
&	\Bigl(I - i\frac{p\tau }{2} \mathfrak{D}_H - i\frac{q \tau}2 \mathfrak{F}^{n+\frac{1}2}\Bigr) \texttt{vec}(U^{n+\frac{1}2}_{i,j}) = \texttt{vec} \bigl( U^n_{i,j} + i\frac{q\tau}2 \Gamma(x_i-y_j) \Phi^{n+\frac{1}2}_{i,j} \bigr), \label{eq:24} \\
&	U^{n+1}_{i,j} = 2U^{n+\frac{1}2}_{i,j} - U^n_{i,j}. \label{eq:25}
\end{align}

The matrices involved are sparse, a property that should be fully exploited in the solution of the linear systems. A spectral method could be used to implement the hyperbolic Laplacian, but this would make the matrix dense, and in many cases would lead to much slower implementation of each timestep, even if fewer points are needed.

It should also be noted that $\Gamma(x_i-y_j)$ should be properly periodized (including near the corners of the domain $[-\frac{L}2,\frac{L}2]^2$).
One way to do this, for  an autocorrelation $\Gamma$ whose effective support (e.g. up to machine error) is well contained in $[-\frac{L}2,\frac{L}2],$ is by using 
\begin{equation}\label{eq:periodgamma765}
\Gamma\Big(   \left( (x-y+\frac{L}2) \mod L \right) \,-\, \frac{L}2\Big)
\end{equation}
instead of $\Gamma(x-y).$ 
This leads to fully periodized $\Gamma.$ This is needed because $\Gamma$ is evaluated outside the primary domain of $[-\frac{L}2,\frac{L}2],$ especially  near  the top-left and bottom-right corners, requiring values of $\Gamma(\xi)$ for $\xi \in [-L,L].$ 

\subsection{Initialization}\label{sec:initimpl}

In order to execute the timestep of \eqref{eq:23}-\eqref{eq:25} for the first time, a value for $\Phi^{-\frac{1}2}_{i,j}$ is needed. In \cite{besse2004relaxation} a naive initialization approach is discussed, which adapted for this problem would amount to the following
\begin{equation}
\mbox{Naive initialization: } \quad \Phi^{-\frac{1}2}_{i,j} =U^0_{i,i}-U^0_{j,j}.
\end{equation}
As we will see, this naive initialization yields the expected second order of convergence for $u,$ but only first order convergence in time for $\phi.$

In \cite{Athanassoulis2023a,besse2021energy}, ideas amounting to a $2^{nd}$ order approximation of $\Phi(-\tau/2)$ are discussed. In this paper we define
\begin{equation}
\begin{array}{c}
	
\left\{ 
\begin{array}{l}
\Phi^{*}_{i,j} =  U^0_{i,i} - U^0_{j,j}, \\
\Bigl(I + i\frac{p\tau }{4} \mathfrak{D}_H + i\frac{q \tau}4 \mathfrak{F}^{*}\Bigr) \texttt{vec}(U^{-\frac{1}4}_{i,j}) = \texttt{vec} \bigl( U^0_{i,j} - i\frac{q\tau}4 \Gamma(x_i-y_j) \Phi^{*}_{i,j} \bigr), \\
U^{-\frac{1}2}_{i,j} = 2U^{-\frac{1}4}_{i,j} - U^0_{i,j},
\end{array}
\right\} \\[18pt]

\mbox{Advanced initialization: } \quad  \Phi^{-\frac{1}2}_{i,j} = U^{-\frac{1}2}_{i,i} - U^{-\frac{1}2}_{j,j}.

\end{array}
\end{equation}

This advanced initialization leads to schemes of the expected order in both space and time, for both $u$ and $\phi.$ 

The question of initialization for relaxation Crank-Nicolson schemes has recently attracted focused attention. This kind of scheme was introduced in the cubic nonlinear Schr\"odinger equation, and in that context there is no loss of order with the naive initialization (not even in the auxiliary variable)  \cite{besse2004relaxation, zouraris2023error}. On the other hand, for other problems loss of order in the auxiliary variable has been reported with the naive initialization \cite{athanassoulis2024efficient,Athanassoulis2023a,zouraris2021error}. Currently,  advanced second order initializations are considered state-of-the-art, even for the nonlinear Schr\"odinger equation \cite{ besse2021energy}.

 To illustrate the numerical behavior for this problem, we use an exact solution and examine the experimental order of convergence (EOC) in space and time with each initialization. Details for the exact solution and precise definition of the errors can be found in Section \ref{sec:EOC}.

\begin{table}
\caption{Time EOC for $u$ and $\phi$  with advanced initialization. The timestep $\tau$ is refined by a factor of $\sqrt{2}$ in each row, while the spatial mesh size $h$ remains fixed. The final time is $T=0.6.$ The errors $\mathcal{E}_u,$ $\mathcal{E}_\phi$ are defined in \eqref{eq:err1}, \eqref{eq:err2}, and  the invariant errors $\delta I_j$  in \eqref{eq:defdeltainv}.}\vspace{-2mm}
\begin{center}
\begin{tabular}{||c|c||c|c||c|c|c|c|c||} \hline
$h$ & $\tau$ & $\mathcal{E}_u$  & EOC  $u$ &  $\mathcal{E}_\phi$  & EOC  $\phi$ & $\delta I_0 + \delta I_1$ & $\delta I_2$ & $\delta I_3$ \\ \hline 
  0.04   &     0.03   &     0.0046064    &          -- &            0.01075            &    --  & $7e-16$ & $5e-7$ & $3e-4$ \\ \hline   
    0.04  &  0.021213   &     0.0023352   &      1.9602  &         0.0053544     &      2.0112   & $1e-16$ & $2e-6$ & $1e-4$ \\ \hline   
    0.04   &    0.015    &    0.0011427    &     2.0621   &        0.0026856      &      1.991   & $9e-16$ & $3e-7$ & $7e-5$ \\ \hline   
    0.04  &  0.010607    &   0.00057461     &    1.9837    &       0.0013398       &    2.0064   & $4e-16$ & $2e-7$ & $4e-5$ \\ \hline   
     \end{tabular}
\end{center}
\label{eoch5}
\end{table}%

\begin{table}
\caption{Space EOC for $u$ and $\phi$ with advanced initialization. The spatial mesh size $h$ is refined by a factor of $2^{1/4}$ in each row, while the temporal step $\tau$ remains constant. The final time is $T=0.6.$ The errors $\mathcal{E}_u,$ $\mathcal{E}_\phi$ are defined in \eqref{eq:err1}, \eqref{eq:err2}, and the invariant errors $\delta I_j$ in \eqref{eq:defdeltainv}. }\vspace{-2mm}
\begin{center}
\begin{tabular}{||c|c||c|c||c|c|c|c|c||} \hline
$h$ & $\tau$ & $\mathcal{E}_u$  & EOC  $u$ &  $\mathcal{E}_\phi$  & EOC  $\phi$ & $\delta I_0 + \delta I_1$ & $\delta I_2$ & $\delta I_3$ \\ \hline 
       0.4   & 0.0005     &    0.013227    &          --    &        0.023294     &         --   & $4e-15$ & $4e-5$ & $8e-4$   \\ \hline   
    0.33636  &  0.0005    &     0.006746   &      3.8854    &        0.011937     &      3.8581  & $4e-15$ & $1e-5$ & $4e-4$    \\ \hline  
    0.28284  &  0.0005    &    0.0033308   &      4.0727    &       0.0059017     &      4.0648  & $1e-15$ & $2e-6$ &  $2e-4$   \\ \hline  
    0.23784  &  0.0005    &    0.0016545   &      4.0379    &       0.0029336     &      4.0338  & $5e-15$ & $6e-7$ & $1e-4$   \\ \hline  
     \end{tabular}
\end{center}
\label{eoch}
\end{table}%

In Tables \ref{eoch5} and \ref{eoch} the expected order in space and time, respectively, is observed when the advanced initialization is used. Table \ref{eoch6} indicates that $\phi$ is clearly first order in time with the naive initialization. Table \ref{eoch2} shows that, with the naive initialization, the spatial order for $\phi$ is obscured in practice by the relatively large time error.  

\subsection{Invariants at the fully discrete level} \label{sec:inarinants}

\begin{proposition}\label{prop:discrinv123}
	Denote
	\begin{equation}\label{eq:disinbvrfty}
	I_0[U^n_{i,j}] = h^2 \sum_{i=0}^{N-1} \sum_{j=0}^{N-1} |\Gamma(x_i-y_j) + U^n_{i,j}|^2  \qquad 
	I_1[U^n_{i,j}] = h \sum_{i=0}^{N-1} U^n_{i,i},
	\end{equation}
	where $U^n_{i,j}$ is defined in \eqref{eq:23}-\eqref{eq:25}. 	Then
	\[
		I_0[U^{n+1}_{i,j}] = 	I_0[U^n_{i,j}], \qquad 
		I_1[U^{n+1}_{i,j}] = 	I_1[U^n_{i,j}].
	\]
\end{proposition}

The proof follows in the same logic as Proposition \ref{prop:bddns}, using summation by parts instead of integration by parts. The definition of the discrete invariant in \eqref{eq:disinbvrfty} essentially corresponds to a composite trapezoid quadrature formula  for periodic functions applied to  \eqref{eq:i0}, 
\eqref{eq:i1}. While the approximation of the continuous invariant values is relatively low order, this value is now propagated in time up to rounding errors (there are large cancellations in the summation by parts so these errors are barely noticeable in some cases). If we opted for a higher order quadrature rule, the approximation of the continuous invariant would be higher order in $h$, but the value wouldn't be preserved in time, as summation by parts does not apply for general quadrature rules.

We use the same quadrature rule for the discretizations of invariants $I_2,$ $I_3$ as well.


\subsection{Experimental order of convergence} \label{sec:EOC}

\begin{figure}
	\includegraphics[width=0.9\textwidth]{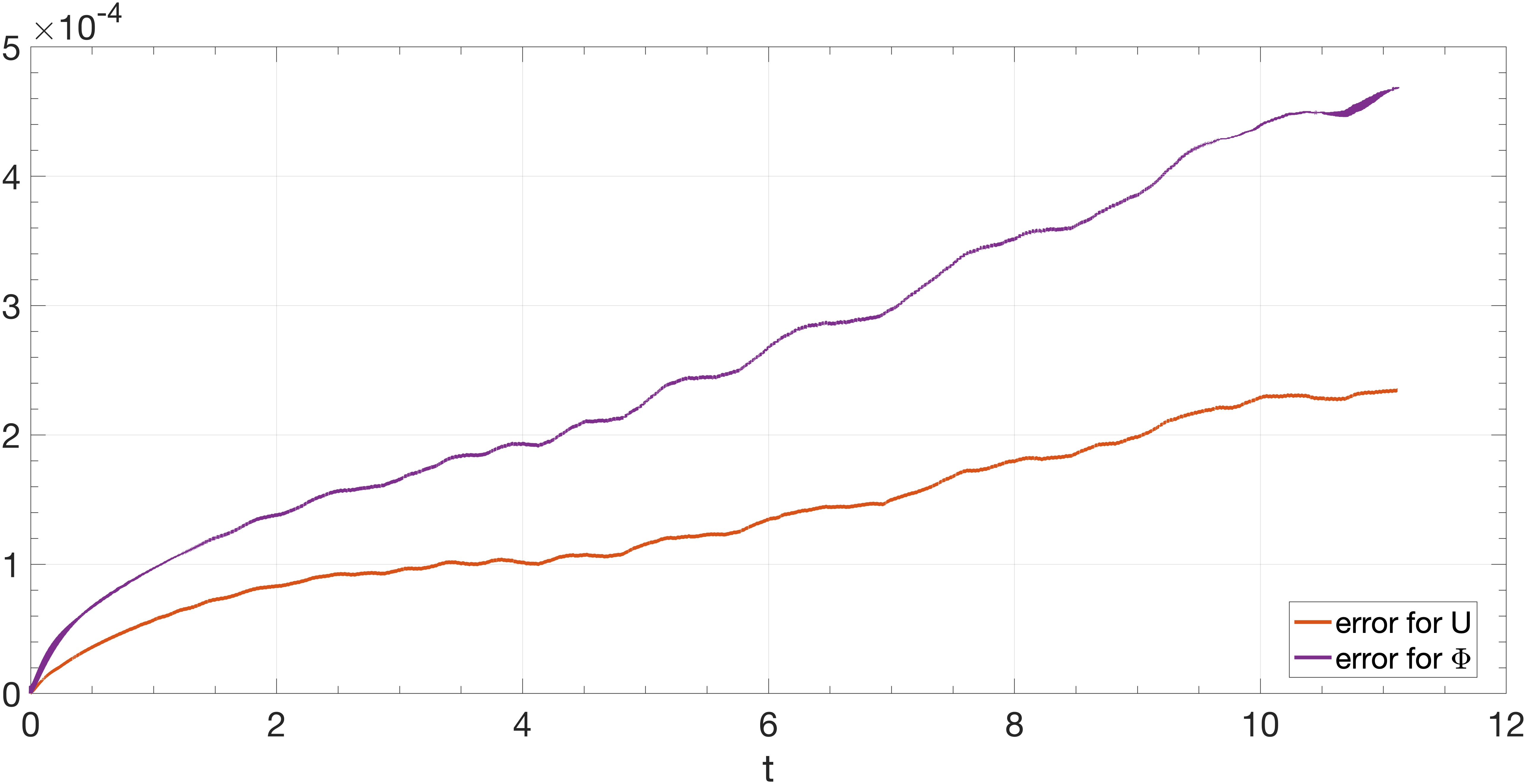}
	\caption{Errors for $u,$ $\phi,$ namely $\mathcal{E}_u^n$ and  $\mathcal{E}_\phi^n$ respectively, plotted against $t^n\in[0,T],$  for the exact solution \eqref{eq:periodizedsoliton}, when the advanced initialization is used. The final time $T=11.1149$ used here amounts to a full lap of the computational domain of length $L=34.4562$, i.e. the soliton returns near the starting position. For this computation the mesh sizes used are $\tau=0.001$ and $h=0.09.$ }\label{Fig:1}
\end{figure}

In the case of zero background, $\Gamma=0,$  exact solutions can be constructed from exact solutions of the underlying nonlinear Schr\"odinger equation. Let $\Gamma=0$ and set  $p=1.7,$ $q=1.1$ in equation \eqref{eq:inhomalb1}. For a soliton amplitude $A\in \mathbb{R}$ and a soliton group velocity $v\in \mathbb{R},$ set $k={v}/{(2p)}$  and $B=A\sqrt{{q}/{(2p)}}.$ Then, it follows that 
\begin{equation}
\begin{array}{c}
u(x,y,t) = A^2 \mathrm{sech}(B(x-vt)) \mathrm{sech}(B(y-vt)) \, e^{ik(x - y)} 
\end{array}
\end{equation}
is an exact solution of \eqref{eq:inhomalb1}. Here we set  $A=1.3,$ $v=3.1.$  For these choices of parameters, $u(x,y,0)$ is numerically zero outside $\Omega := [-\frac{L}2,\frac{L}2]^2$ for $L=10 \pi /k\approx 34.45.$ Thus, when periodic boundary conditions are used on $\Omega,$ the exact solution can be periodized to
\begin{equation}\label{eq:periodizedsoliton}
\begin{array}{c}
w(\xi) =   \left( (\xi+\frac{L}2) \mod L \right) \,-\, \frac{L}2 \\[4pt]
u(x,y,t) = A^2 \mathrm{sech}\bigl( B w(x-vt) \bigr) \mathrm{sech}\bigl(B w (y-vt)\bigr) \, e^{ik(x - y)} 
\end{array}
\end{equation}
 potentially allowing very long times to be considered. This is the exact solution on which  Tables \ref{eoch5} through \ref{eoch2} are based, along with Figures \ref{Fig:1} and \ref{Fig:2}. Observe that the solution has modulus of $O(1)$  in value and derivatives, and this stays so for all times, making the interpretation of the errors straightforward.

The experimental order of convergence (EOC) values in row $i$ of each Table are computed from the runs in rows $i$ and $i-1.$ The errors are defined as
\begin{align}
& \mathcal{E}_u := \max\limits_n \max\limits_{i,j} \left|U^n_{i,j} - u(x_i,y_j,t^n)\right| \label{eq:err1} \\
& \mathcal{E}_\phi := \max\limits_n  \max\limits_{i,j} \left|\Phi^{n-\frac{1}2}_{i,j} - \phi(x_i,y_j,t^{n-\frac{1}2})\right| \label{eq:err2}
\end{align}
with respect to the exact solution  \eqref{eq:periodizedsoliton}. Note that 
 $t^N=T\approx 0.6$ is the final time. Moreover, the relative errors in the preservation of the invariants  are defined as
 \begin{equation}\label{eq:defdeltainv}
 	\delta I_j := \left| \frac{I_j[U^N] - I_j[u_0]}{I_j[u_0]} \right|,
 \end{equation}
 where the trapezoid quadrature rule is used (similarly to Proposition \ref{prop:discrinv123}).

\begin{table}
\caption{Time EOC for $u$ and $\phi$ with naive initialization.  The timestep $\tau$ is refined by a factor of $\sqrt{2}$ in each row, while the spatial mesh size remains constant. The final time is $T=0.6.$ The errors $\mathcal{E}_u,$ $\mathcal{E}_\phi$ are defined in \eqref{eq:err1}, \eqref{eq:err2} and the invariant errors $\delta I_j$ in \eqref{eq:defdeltainv}. }\vspace{-2mm}
\begin{center}
\begin{tabular}{||c|c||c|c||c|c|c|c|c||} \hline
$h$ & $\tau$ & $\mathcal{E}_u$  & EOC  $u$ &  $\mathcal{E}_\phi$  & EOC  $\phi$ & $\delta I_0 + \delta I_1$ & $\delta I_2$ & $\delta I_3$ \\ \hline 
 0.04     &   0.03      &   0.005322       &       --        &    0.09592       &          0    &  $9e-16$  &   $5e-4$   & $7e-4$ \\ \hline 
    0.04  &  0.021213   &     0.0026947    &     1.9637      &     0.066486     &       1.0576  &   $3e-16$ &  $1e-4$    & $1e-4$ \\ \hline  
    0.04  &     0.015   &     0.0013247    &     2.0489      &     0.046349     &        1.041   &  $9e-16$  &  $1e-4$    & $1e-4$ \\ \hline 
    0.04  &  0.010607   &    0.00066544    &     1.9866      &     0.032443     &       1.0293   &  $1e-15$  &  $4e-5$    & $4e-5$  \\ \hline 
     \end{tabular}
\end{center}
\label{eoch6}
\end{table}%

\begin{table}
\caption{Space EOC for $u$ and $\phi$ with naive initialization.  The spatial mesh size $h$ is refined by a factor of $2^{1/4}$ in each row, while the temporal step $\tau$ remains constant. The final time is $T=0.6.$ The errors $\mathcal{E}_u,$ $\mathcal{E}_\phi$ are defined in \eqref{eq:err1}, \eqref{eq:err2} and the invariant errors $\delta I_j$ in \eqref{eq:defdeltainv}. }\vspace{-2mm}
\begin{center}
\begin{tabular}{||c|c||c|c||c|c|c|c|c||} \hline
$h$ & $\tau$ & $\mathcal{E}_u$  & EOC  $u$ &  $\mathcal{E}_\phi$  & EOC  $\phi$ & $\delta I_0 + \delta I_1$ & $\delta I_2$ & $\delta I_3$ \\ \hline 
       0.4   & 0.0005    &     0.013227     &         --      &      0.023718     &           --  & $4e-15$ & $4e-5$  &  $8e-4$   \\ \hline 
    0.33636  &  0.0005   &     0.0067462    &     3.8854      &       0.01231     &      3.7846   & $3e-15$ & $1e-5$  &  $4e-4$   \\ \hline 
    0.28284  &  0.0005   &     0.0033309    &     4.0726      &     0.0062821     &      3.8821   &  $8e-15$ & $2e-6$ &  $2e-4$   \\ \hline 
    0.23784  &  0.0005   &     0.0016547    &     4.0375      &     0.0033981     &       3.546   & $3e-15$ & $7e-7$ &  $1e-4$   \\ \hline 
     \end{tabular}
\end{center}
\label{eoch2}
\end{table}%

The coarsest $h$ used is relatively large at $h\approx 0.4,$ but even so we get  two significant digits correct in both $u$ and $\phi$ in all cases. 


\subsection{Qualitative behavior of the error}

\begin{figure}
	\includegraphics[width=0.9\textwidth]{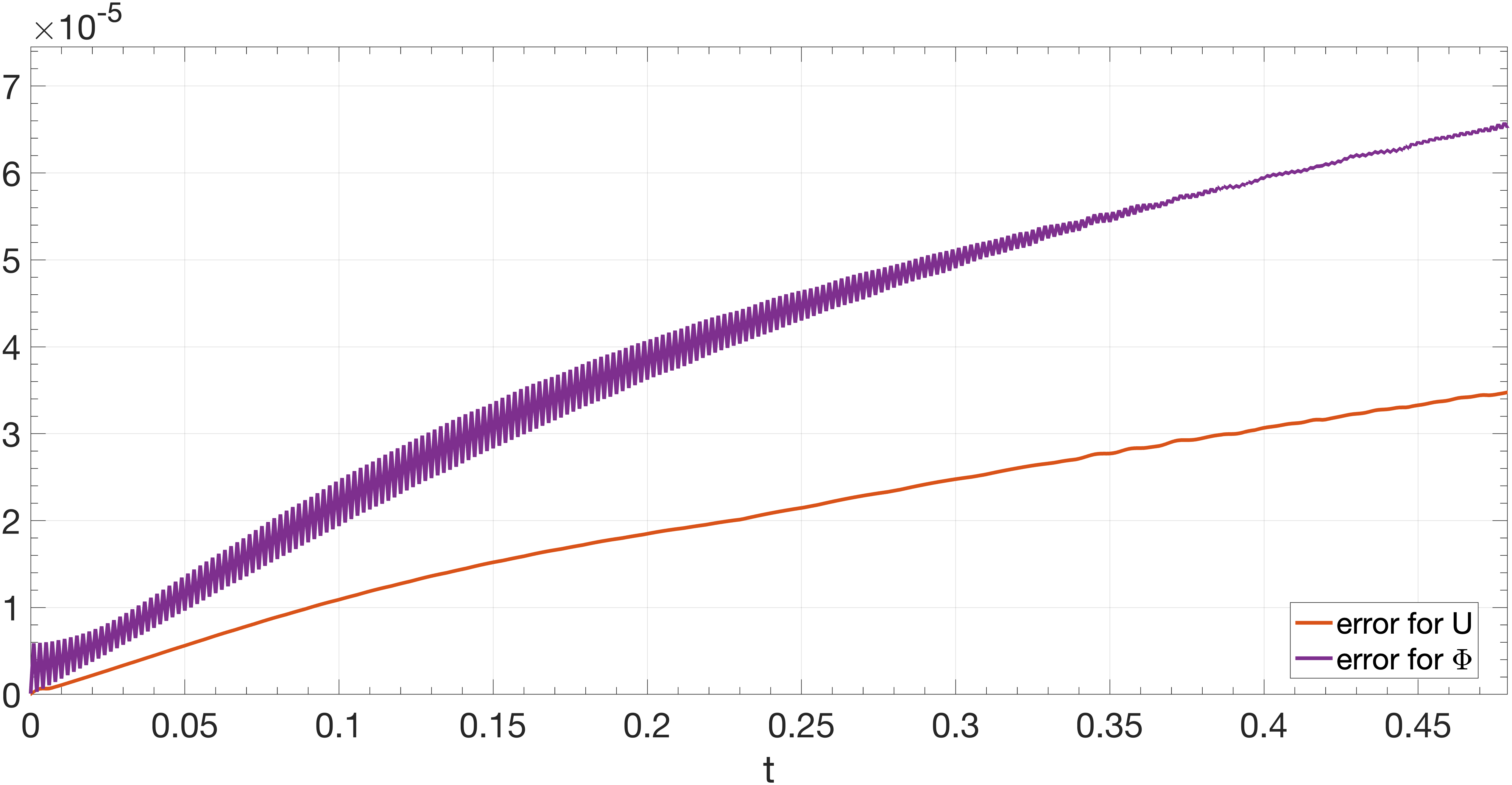}
	\caption{Zoom in of the error in Figure \ref{Fig:1} at the initial stage. The oscillations relax after some time. The naive initialization produces substantially more pronounced oscillations. }\label{Fig:2}
\end{figure}

\begin{figure}
	\includegraphics[width=0.9\textwidth]{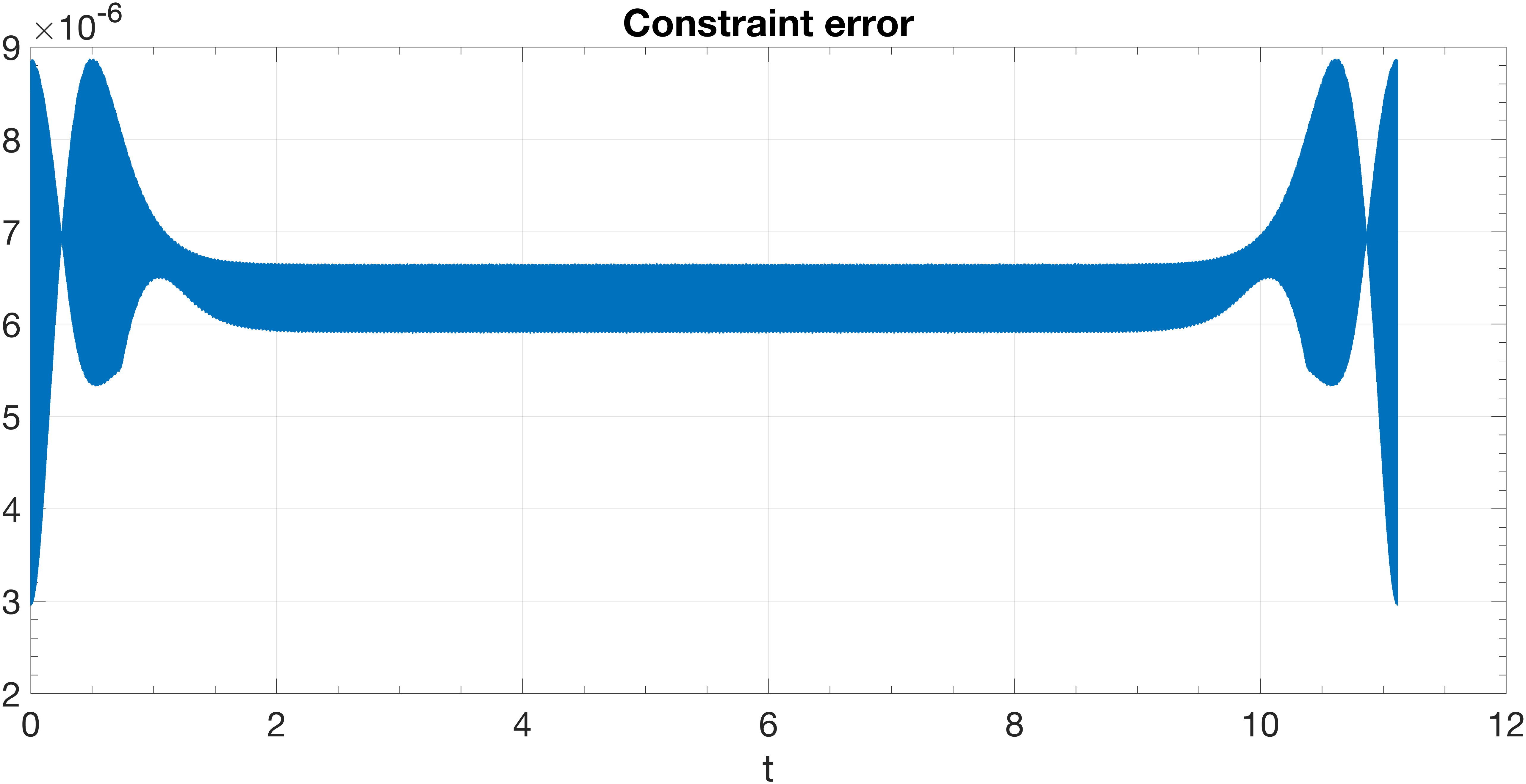}
	\caption{Constraint error, as defined in equation \eqref{eq:constrerrdef}, for the example of Figure \ref{Fig:1}. }\label{Fig:3}
\end{figure}

Some more detailed observations can be made for the qualitative behavior of the error on a longer-time run, $T\approx11.1.$ Here we use again the soliton-type solution of \eqref{eq:periodizedsoliton}, with the advanced initialization and well-resolved discretization, cf.  Figure \ref{Fig:1}. After some time, the soliton exits from the top right corner of the computational domain, and re-enters from the bottom left corner, until it returns to the origin. The timestep is $\tau=0.001$ and the spatial mesh size $h=0.09.$ 

It is seen that the error in $U$ is not very sensitive to the error in $\Phi:$ the latter grows faster and to larger values, cf. Figure \ref{Fig:1}. 

Another interesting phenomenon is that, after a rapid initial growth, the error seems to grow roughly linearly in time (and not exponentially, despite the problem being nonlinear). 

If we look at more fine behaviours of the error, we observe  oscillations of the exact error early on, cf. Figure \ref{Fig:2}. Whenever the naive initialization is used instead of the advanced one, these oscillations become more pronounced.

This scheme rests crucially on the introduction of an auxiliary variable on a staggered time-grid, $\Phi^{n+\frac{1}2}$ which aims to represent $u(x,x,t^{n+\frac{1}2}) - u(y,y,t^{n+\frac{1}2}).$ The direct relationship between $\Phi^{n+\frac{1}2}$ and $U^n$ is through the first equation in \eqref{eq:rCN1},  i.e., the $U^n$ interpolates appropriately in time the $\Phi^{n\pm\frac{1}2}.$ To see how well the relationship between the main and auxiliary variables holds up, one may ask how well the reverse holds, i.e. how well $\Phi^{n+\frac{1}2}$ interpolates $U^n,$ $U^{n+1}.$ This is measured by the quantity
\begin{equation}\label{eq:constrerrdef}
	\mathcal{E}_{constr}^n := \max\limits_{i,j} \left| \frac{1}2 \left( U^{n+1}_{i,i}-U^{n+1}_{j,j} - U^n_{i,i} + U^n_{j,j}   \right) - \Phi^{n+\frac{1}2}_{i,j} \right|.
\end{equation}
To the best of our knowledge, such measures of consistency for auxiliary variables do not have a standard name, and we call this the ``constraint error''. This seems to fluctuate initially (roughly over the same time period as the oscillations of Figure \ref{Fig:2}), and then settles to a constant value, much smaller than the error -- see Figure \ref{Fig:3}. Near the final time, when the soliton returns near its original position, we see a ``refocusing'' of the numerical artifacts -- most clearly in Figure \ref{Fig:3}, but also some oscillations in the error in $\Phi$ in Figure \ref{Fig:1}. This is a consequence of the periodization.

The apparent boundedness of this constraint error in time is a promising feature for the use of the scheme in long times, and contrasts  other ``auxiliary variable'' types of schemes  where analogous quantities are known to grow \cite{Alsafri2023}.

\begin{figure}
	\includegraphics[width=0.9\textwidth]{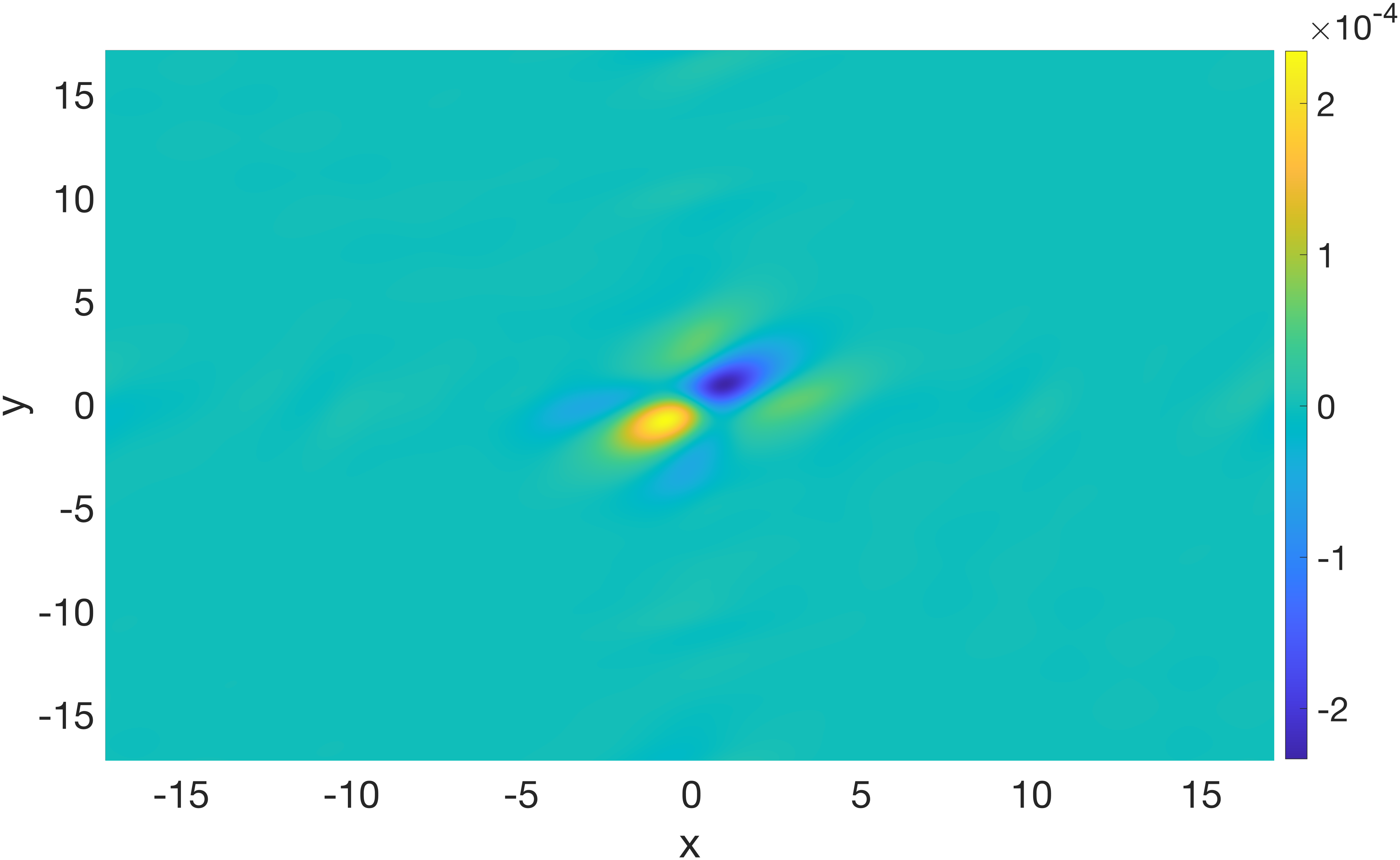} 
	\caption{This plot is the real part of $u(t^N)-U^N$ over $(x,y)\in[-\frac{L}2,\frac{L}2]^2$ at final time $T\approx 11.1,$ for the run described in Figure \ref{Fig:1}. The shape of the soliton is very well preserved; the bulk of the error is due to the phase of the numerical solution staying slightly behind that of the exact solution.}
	\label{Fig:4}
\end{figure}

Finally, by examining the full difference between the numerical and exact solutions, we can get a sense of where the bulk of the error comes from.  In Figure \ref{Fig:4} we see that the main part of the error in this problem is in fact due to the phase of the numerical solution staying slightly behind the exact solution, while the shape and location of the soliton are reproduced to very high accuracy. 

\medskip 

Invariants $I_0, I_1$ are preserved to machine error, and invariants $I_2, I_3$ are approximately preserved to high accuracy in the example of Figure \ref{Fig:1}; see Table \ref{Tab:5}. 

\begin{table}
\caption{Invariant errors, as defined in equation \eqref{eq:defdeltainv}, for the example of Figure \ref{Fig:1} ($T\approx 11.1,$ $\tau=0.001$ and $h=0.09$). }\vspace{-2mm}
\begin{center}
\begin{tabular}{||c|c|c|c||} \hline
$\delta I_0$ & $\delta I_1$ & $\delta I_2$ & $\delta I_3$  \\ \hline 
  $3e-16$ & $3e-15$ & $4e-10$ & $5e-6$     \\ \hline 
\end{tabular}
\end{center}
\label{Tab:5}
\end{table}%

\section{Landau damping and modulation instability} \label{sec:LD}

\subsection{The stability condition}\label{subsec:41}

A central feature of the Alber equation \eqref{eq:inhomalb1} is the bifurcation between a linearly stable and a linearly unstable regime. The stable regime is similar to Landau damping \cite{Athanassoulis2018}, and the inhomogeneity $u$ disperses despite the presence of the homogeneous background $\Gamma$ which would be expected to act as a source term. In the unstable case, the inhomogeneity $u$ triggers unstable modes that grow exponentially -- at least initially \cite{Alber1978}. This  is often called modulation instability, in analogy with the standard modulation instability of the nonlinear Schr\"odinger equation (NLS) \cite{benjamin1967disintegration,zakharov1968stability}.

The linear  instability condition is formulated in terms of the background power spectrum of the problem, $P(k) = \mathcal{F}_{x\to k} [\Gamma(x)]$  \cite{Alber1978}, and it takes the form of a Penrose-type condition. The sufficient condition for instability  is \cite{Alber1978,Athanassoulis2018}
\begin{equation}\label{eq:instability}
	\exists X_* \qquad \exists \omega_* \mbox{ with } \Rea(\omega_*)>0 \qquad \mathbb{H}[D_{X_*}P](\omega_*) = 4\pi\frac{p}q,
\end{equation}
where $\mathbb{H}$ is the Hilbert transform and
 $D_XP(k)$ is the divided difference of the background spectrum $P$ with increment $X,$
\begin{equation}\label{eq:DXP}
	D_X P(k)=\left\{\begin{array}{cc}
\frac{P\left(k+\frac{X}{2}\right)-P\left(k-\frac{X}{2}\right)}{X}, & X \neq 0, \\[6pt]
P^{\prime}(k), & X=0.
\end{array}\right.
\end{equation}

\begin{figure}
	\includegraphics[width=0.9\textwidth]{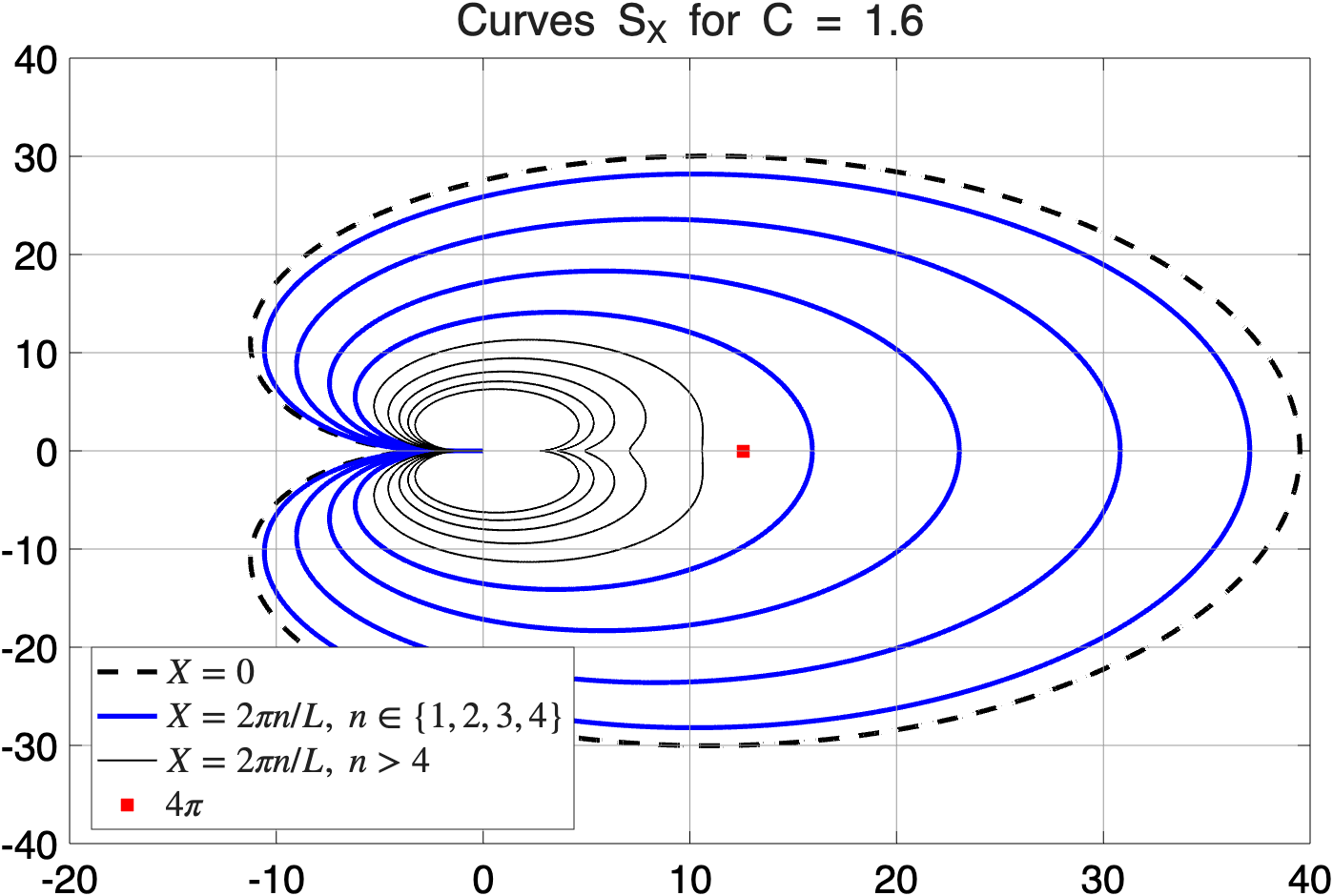} 
	\caption{The curves $S_X(t),$ defined in \eqref{eq:curveS}, are plotted on the complex plane for various values of $X\in\mathbb{R}$ as in \eqref{eq:curveS}. In the context of Section \ref{subsec:42}  $p=q=1,$ hence the target point is $4\pi p /q =4\pi.$ For this particular case $C=1.6,$ the first four harmonics of the computational domain, $X=2\pi n / L,$ $n\in\{1,2,3,4\},$ are unstable. The bandwidth of unstable wavenumbers for this case is approximately $1.13.$}
		\label{Fig:Nyquist}
\end{figure}

The  instability condition \eqref{eq:instability} is not straightforward to check, as it involves the existence or not of  solutions in a system of two equations (real and imaginary parts of $\mathbb{H}[D_XP](\omega)=4\pi p/q$) in three unknowns ($X,\Rea(\omega),\Imag(\omega)$). A constructive way to check stability was introduced in \cite{Athanassoulis2018}: the curve
\begin{equation}\label{eq:curveS}
	S_X(t) := \mathbb{H}[D_XP](t) - i D_XP(t), \qquad t \in \mathbb{R},
\end{equation}
is drawn on the complex plane. This is a closed curve starting and ending at $0$ (and it might be self-intersecting for multimodal spectra $P$). If the winding number of $S_X$ around $4\pi p/q$ is nonzero, then the wavenumber $X$ is unstable. This Nyquist-type  method computes the solutions $X_*$ of \eqref{eq:instability}, but it only clarifies the existence of $\omega_*$ -- not their values. In Figure \ref{Fig:Nyquist} this technique is demonstrated for the Gaussian spectrum \eqref{eq:gamma} with $C=1.6.$ From this, the bandwidth of unstable wavenumbers can be computed. For the required computation of the Hilbert transform we use Weideman's rational eigenfunction method \cite{weideman1995computing}. Standalone scripts implementing this check of the stability condition can be found in \url{https://github.com/aathanas/penrose4alber}.

If condition \eqref{eq:instability} holds for some $X_*,$ then that $X_*$ is an unstable wavenumber. The instability  is a baseband modulation instability in the sense that  unstable wavenumbers (when they exist) are of the form $|X_*|< X_{\text{max}},$ i.e., large enough wavenumbers are always stable. Moreover, the bandwidth of unstable wavenumbers, $2X_{\text{max}},$ controls the lengthscales in the emerging coherent structure (larger bandwidth leads to more localized coherent structures), while the largest $\Rea(\omega_*)$ controls the rate of growth of the inhomogeneity, according to the linear stability analysis \cite{Athanassoulis2017RogueWaves}. An important question is whether the fully nonlinear dynamics remain close to the linearized dynamics, especially for physically realistic inhomogeneities.  For example, many nonlinear Landau damping results have extremely stringent smallness conditions for the inhomogeneity and its derivatives up to high order, which are questionable in any real-world system.

\medskip 

A question that is debated in the ocean waves community is the nature of the transition from stable to unstable behavior, and especially how  it manifests phenomenologically in the ocean waves context. To fix ideas, consider the one-parameter family of  spectra of the form $P(k)= C P_0(k)$ for some smooth, unimodal profile $P_0(k).$  It is easy to see that there exists a critical value $C_*$ such that if $C<C_*$ the spectrum is stable (in the sense of condition \eqref{eq:instability}), while for $C>C_*$ it is unstable. The question then becomes, how different are the solutions of \eqref{eq:inhomalb1} for $C=C_*-\epsilon$ versus $C=C_*+\epsilon$ for some $0<\epsilon\ll 1?$ Is there a striking change in the qualitative behavior of solutions around $C=C_*?$ 

A key piece of context  for this question is the timescales of interest: physically, we expect the original equation \eqref{eq:inhomalb1} to be valid for an $O(1)$ timescale $T,$ but not for $t\to\infty.$ In that context, instability that takes an asymptotically long time to manifest would appear similar to stability when $t\in[0,T]$ \cite{Athanassoulis2017RogueWaves}.  In Section \ref{subsec:42} we investigate  numerically  the growth in the case of weak instability. 

Another aspect of this question is how large the inhomogeneity eventually becomes. In the ocean waves context, $u$ is compared to an $O(1)$ homogeneous background. Thus, even if $u$ itself changes a lot, as long as it remains dominated by the background then the phenomenology of the sea state does not change noticeably. We develop an idea of \cite{stiassnie2008recurrent} and introduce the {\em inhomogeneity amplification factor} (IAF; how much the inhomogeneity grows compared to its initial condition) and the {\em total amplification factor} (TAF; how large are the localized events of the total sea state compared to the homogeneous background). Modulation instability leads to large IAF, but if the TAF is small, then the modulation instability does not produce phenomenologically relevant effects.  This is discussed in some detail in Section \ref{sec:4.3}.



A subtle issue in this discussion is the localization of the problem from the full space $(x,y)\in \mathbb{R}^2$ to a finite domain $(x,y)\in[-\frac{L}2,\frac{L}2]^2.$ The linear stability analysis is typically carried out on the full space. In   \cite{Athanassoulis2023b} the question of how to localize the Alber equation and its linear stability analysis on a finite computational domain $(x,y)\in[-\frac{L}2,\frac{L}2]^2$ was investigated.  It was found that setting
\begin{equation}\label{eq:gammasoec}
	\Gamma(x) = \frac{1}L \sum\limits_{n\in\mathbb{Z}} P(\frac{n}L) e^{2\pi i \frac{n \cdot x}L}
\end{equation}
leads to the same stability condition as full space, provided that $L$ is large enough.  In other words, we start from the continuous power spectrum $P(k)$ (often measured in the field in the ocean waves context), and  generate an autocorrelation function $\Gamma$  according to \eqref{eq:gammasoec}. This in turn provides a benchmark for the size of $L$ required:  the effective support of $\Gamma_{\text{inf}}(x):= \mathcal{F}^{-1}_{k\to x}[P(k)]$ must be well contained in $[-\frac{L}2,\frac{L}2].$ As long as that is true,  then equations  \eqref{eq:gammasoec}  and \eqref{eq:periodgamma765} are  numerically interchangeable.

\subsection{Numerical investigation of a model spectrum}\label{subsec:42}

The onset of modulation instability is not yet fully understood, but it is considered one of the possible mechanisms behind the growth of oceanic rogue waves \cite{dysthe2008oceanic}. Hence, investigating its fully nonlinear behavior presents great interest  from a theoretical dynamical systems point of view, as well as an ocean waves point of view. Capturing the growth (or not) of the inhomogeneity is precisely why the discrete balance law of Proposition \ref{prop:discrconbll} is important, as it guarantees that the numerical solution qualitatively follows  the exact one.
Moreover, the second-order in time and fourth-order in space of the numerical method presented here make tractable the large-domain and long-time computations required to conclusively investigate  nonlinear Landau damping (in the stable case), and the growth of coherent structures (in the unstable case). 

For clarity and compactness, here we investigate problems where all the terms are of comparable size, $p=q=1,$ and use a Gaussian power spectrum  
\begin{equation}\label{eq:gamma}
	P(k)  = \frac{C^2}{\sigma}   e^{-\pi \frac{k^2}{\sigma^2}}, \qquad \sigma=0.36, \qquad C \in [0.9,2.2].
\end{equation}
The periodized autocorrelation that corresponds to this spectrum is
\begin{equation}\label{eq:gamma2}
	\Gamma_{\mathrm{inf}}(y)  = C^2   e^{-\pi \sigma^2 y^2 }, \qquad  \Gamma(y) = \Gamma_{\mathrm{inf}}\Big(   \left( (y+\frac{L}2) \mod L \right) \,-\, \frac{L}2\Big)
\end{equation}

The parameter $C$ controls the RMS wave amplitude  in the background sea state. Moreover, this is a relatively  narrow power spectrum, i.e.,  the situation in which modulation instability is expected to appear. More specifically, the stability condition shows onset of instability at $C\approx0.99.$

We use an initial inhomogeneity   of the form
\begin{equation}\label{eq:u0}
\begin{array}{c}
f_0(x,y)=0.05 \cdot   e^{- 0.06x^2 - 0.07y^2}   ( 1+ A_1 \cos(0.3x)\cos(0.2y) +A_2 x + A_3y ), \\[4pt]
u_0(x,y) = \frac{1}2 \left( f_0(x,y) + \overline{f_0(y,x)} \right).
\end{array}
\end{equation}
In order to have directly comparable results for different values of $C,$ we used 
\begin{equation}\label{eq:u01}
A_1=0.3 +  0.8 i, \qquad 
A_2=-0.2, \qquad  A_3=0.1 i
\end{equation}
in all the runs discussed in Section \ref{subsec:42}. Values drawn from a complex normal distribution for $A_1,A_2,A_3$ yield similar results, i.e., the results are not particularly sensitive to the exact values -- as we will see in Section \ref{sec:4.3}, where randomized values are used. The strongest dependence of the results on the initial data is through the overall size of the initial condition, as for smaller initial inhomogeneities it can take substantially longer for the inhomogeneity to grow in the unstable case. The initial datum of \eqref{eq:u0}, \eqref{eq:u01} has $\|u_0\|_{L^2}\approx 0.31$ and $\|u_0\|_{L^\infty_{x,y}}\approx 0.07;$ in other words it is small but not vanishingly so. We use an initial inhomogeneity that is neither  symmetric nor a very simple function to highlight the robustness of the qualitative behavior of the problem. 

\begin{figure}
	\includegraphics[width=0.9\textwidth]{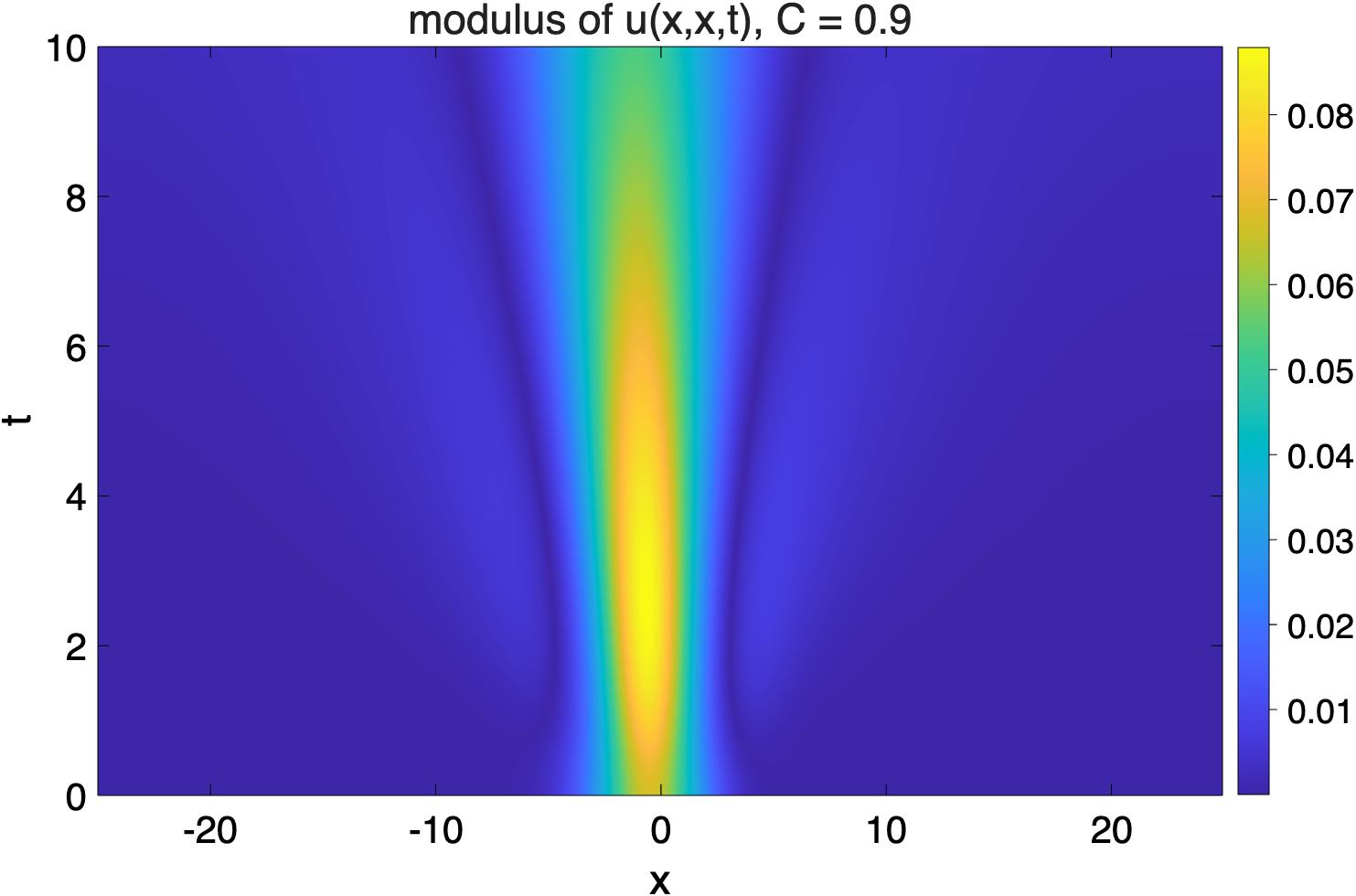}
	\caption{The modulus  $|u(x,x,t)|$ is plotted. The inhomogeneity $u$ satisfies the fully nonlinear equation \eqref{eq:aug1} with $p=q=1,$ $\Gamma$ consistent with equations \eqref{eq:gamma}, \eqref{eq:gamma2} and $C=0.9.$ The initial condition $u_0$ is given in \eqref{eq:u0}, \eqref{eq:u01}. Discretization parameters are given in \eqref{eq:dxdt}. The invariant relative errors for this computation are $\delta I_0 = 1e-15,$ $\delta I_1 = 1.8e-16, \delta I_2 = 1.5e-7$ and $\delta I_3 =1.8e-2.$ In other words $I_0,$ $I_1$ are indeed conserved at the discrete level, and $I_2, I_3$ are reasonably well preserved. It must be noted that $u$ does not grow meaningfully  from its initial condition in either $L^2$ or $L^\infty$ sense.}
		\label{fig:8}
\end{figure}

\begin{figure}
	\includegraphics[width=0.9\textwidth]{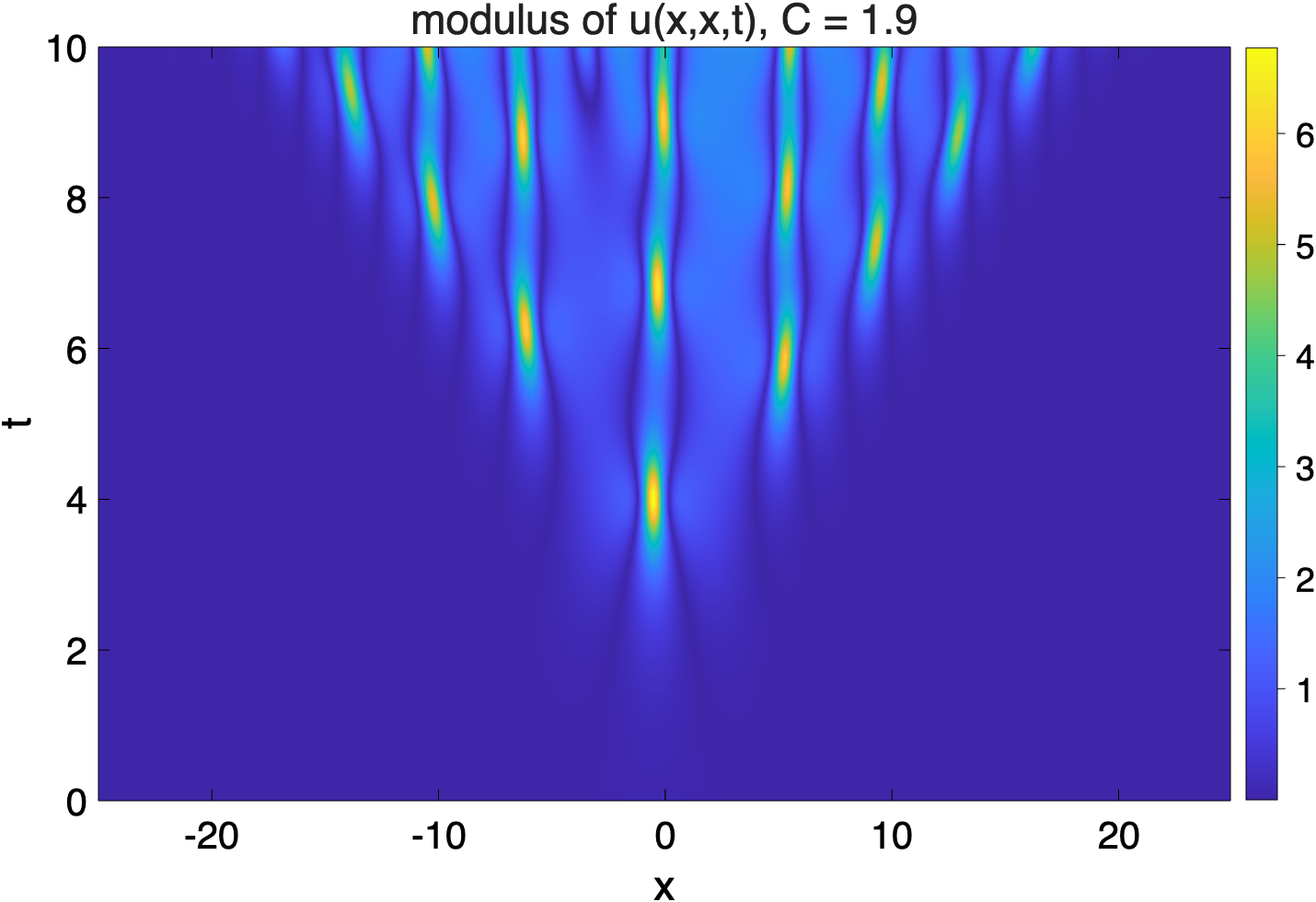}
	\caption{The modulus  $|u(x,x,t)|$ is plotted. The inhomogeneity $u$ satisfies the fully nonlinear equation \eqref{eq:aug1} with $p=q=1,$ $\Gamma$ consistent with equations \eqref{eq:gamma}, \eqref{eq:gamma2} and $C=1.9.$ The initial condition $u_0$ is given in \eqref{eq:u0}, \eqref{eq:u01}. Discretization parameters are given in \eqref{eq:dxdt}. The invariant relative errors for this computation are $\delta I_0 = 1e-15,$ $\delta I_1 = 1e-13, \delta I_2 = 5e-3$ and $\delta I_3 >1.$ In other words $I_0,$ $I_1$ are indeed conserved at the discrete level, $I_2$ is reasonably well preserved, while $I_3,$ which involves second derivatives, is not well preserved. It must be noted that $u$ grows two orders of magnitude from its initial condition in either $L^2$ or $L^\infty$ sense. }
	\label{fig:3}
\end{figure}

All numerical results in this Section, including Figures \ref{fig:3}-\ref{fig:8}, are carried out with $\Gamma$ as in equations \eqref{eq:gamma}, \eqref{eq:gamma2} ($C$ specified in each case), $u_0$ as in \eqref{eq:u0}, \eqref{eq:u01}, with $(x,y)\in[-\frac{L}2,\frac{L}2]^2$ and
\begin{equation}\label{eq:dxdt}
 L=50, \quad \mbox{ timestep } \tau=10^{-3}, \mbox{ spatial mesh size } h=9\cdot10^{-2}.
\end{equation}

\medskip

In terms of visualizing the results, there is no need to plot snapshots of a two-dimensional function. Since the potential $V(x,t)=u(x,x,t)$ controls  the size of  deviation  from the homogeneous background at location $x$ and time $t,$ we plot $|u(x,x,t)|$ against space $x$ and time $t.$ We also plot the growth of certain norms of $u$ in time.

In Figure \ref{fig:8} we see a typical example of nonlinear Landau damping for the equation \eqref{eq:inhomalb1}. There is no meaningful growth of the inhomogeneity $u$ present, only dispersion. 

\medskip

If $C$ increases sufficiently, then we get  modulation instability. In that case, the inhomogeneity first goes through a phase of exponential growth, and a localized coherent structure with values of $O(1)$ appears. Then, growth in $L^\infty$ sense stops, and a distinctive pattern emerges, consisting of similar localized structures spreading outwards in space to form a lattice within a space--time cone. A typical example of the space-time lattice of coherent structures can be seen in Figure \ref{fig:3}. This matches very closely  the universal behavior of modulation instability in other wave problems \cite{biondini2018universal}. The main features of this cone-supported lattice (such as the size and distance of the peaks, the maximum amplitude etc) seem to depend mainly on the homogeneous background $\Gamma,$ and not on the initial condition. This finding is further supported by the systematic Monte Carlo investigation reported in Section \ref{sec:4.3}. To the best of our knowledge, there are no results that identify this meta-stable pattern directly  (i.e. without solving the IVP \eqref{eq:aug1} for  particular initial inhomogeneities).

\begin{figure}
	\includegraphics[width=0.9\textwidth]{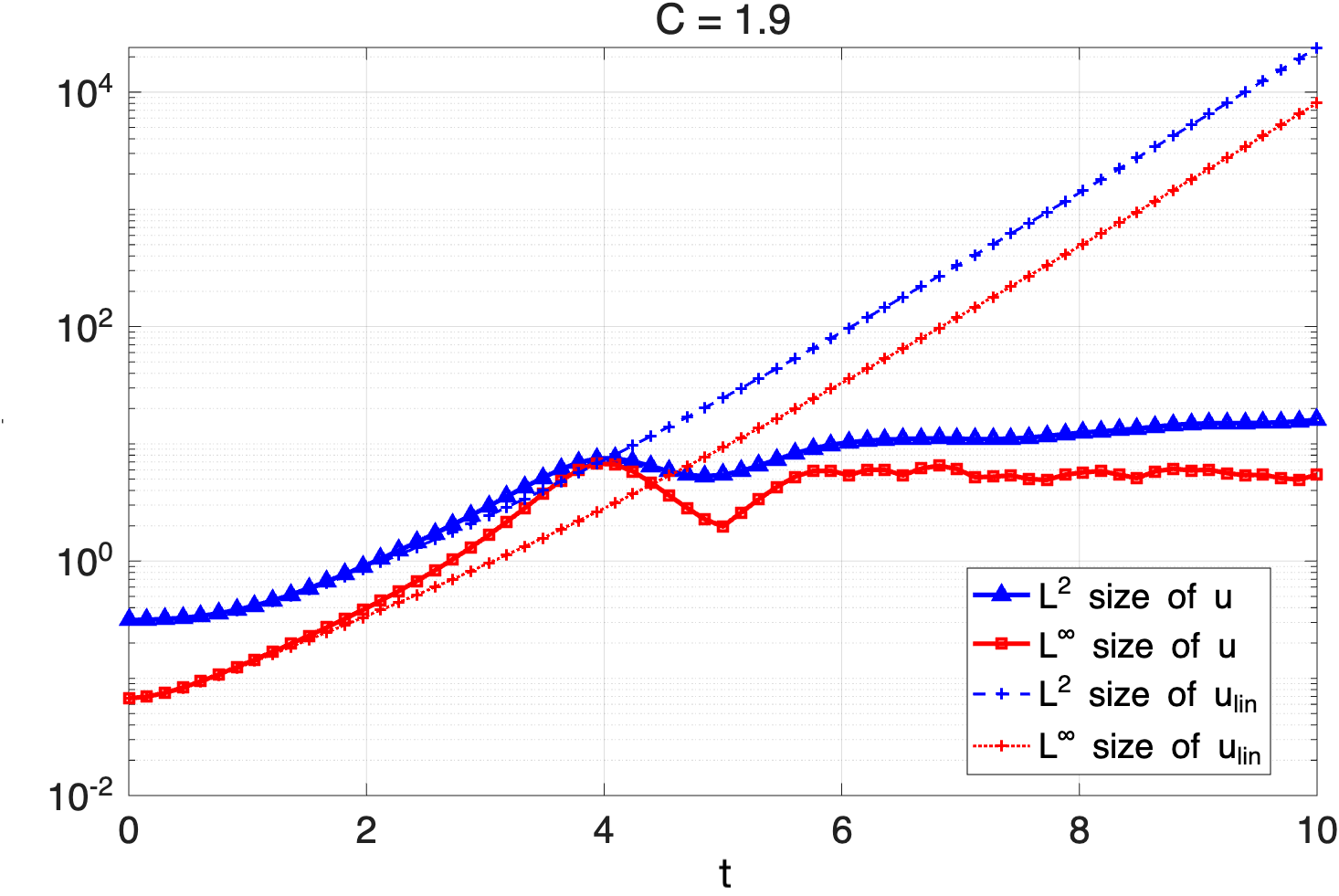}
	\caption{Evolution in time of the norms of the solution $u$ of the fully nonlinear problem \eqref{eq:aug1};  and $u_{\mathrm{lin}},$ the solution of the linearized problem \eqref{eq:inhomalb1LIN}; in a case with modulation instability. All parameters are as in Figure \ref{fig:3}. }
		\label{fig:4}
\end{figure}

When $C=1.9$ as in Figure \ref{fig:3}, this is  a  relatively intense  background  giving rise to a relatively strong instability. Here the values of $|u|$  peak at about $6,$  clearly dominating the background of $1.9.$ The two distinct stages are visible in Figures \ref{fig:3} and \ref{fig:4}: first, there is exponential growth of the inhomogeneity to a single, localized coherent structure.  Once the first coherent structure is formed, there is no more growth in the $L^\infty$ norm, but the pattern propagates outwards forming the lattice-within-a-cone.

We call this a meta-stable pattern because it is robust with regard to small changes in  $p,$ $q$ and $\Gamma,$ and even with large relative changes in the initial condition $u_0:$ any small initial condition we tried  gave rise to a cone like in Figure \ref{fig:3}, with similar maximum value and similar spacing on time-space.

For purposes of comparison with the linear stability analysis of the problem, we also solve the linearized Alber equation,
\begin{equation}\label{eq:inhomalb1LIN}
\begin{array}{c}
		i \partial_t u_{\text{lin}} + p (\Delta_x - \Delta_y) u_{\text{lin}} + q\Bigl(u_{\text{lin}}(x,x,t)-u_{\text{lin}}(y,y,t)\Bigr) \,  \Gamma(x-y)  =0,  \\[8pt]
 \qquad u_{\text{lin}}(x,y,0) = u_0(x,y).
\end{array}
\end{equation}
In Figure \ref{fig:4} we see the evolution in time of the $L^2$ and $L^\infty$ norms for $u$ and $u_{\mathrm{lin}}$ for $C=1.9.$ By $t=1$ or so the solutions of both the fully nonlinear and linearized equations have reached their fastest rate of growth. This picking-up-steam stage can take much longer when the instability is weaker (e.g., smaller $C$). 
Figure \ref{fig:4} shows that the linearized equation \eqref{eq:inhomalb1LIN} captures very well the initial stage of exponential growth of the unstable modes. On the other hand, the linearized equation completely misses the maxing out of the $L^\infty$ norm and the formation of a coherent structure. In conclusion, the linear stability analysis is found to accurately predict the onset of modulation instability, as well as the rate of growth in the early stage. However, it clearly misses features like the maximum amplitude  $\max|u|$, the space-time cone etc. Full solution of the nonlinear equation \eqref{eq:aug1} is required for that.

\begin{figure}
	\includegraphics[width=0.9\textwidth]{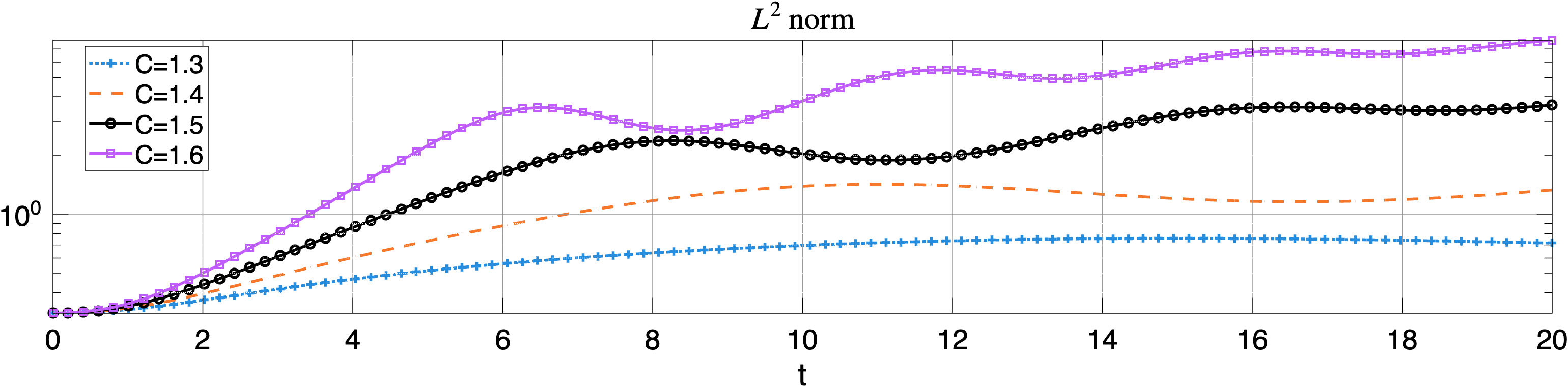}
	\includegraphics[width=0.9\textwidth]{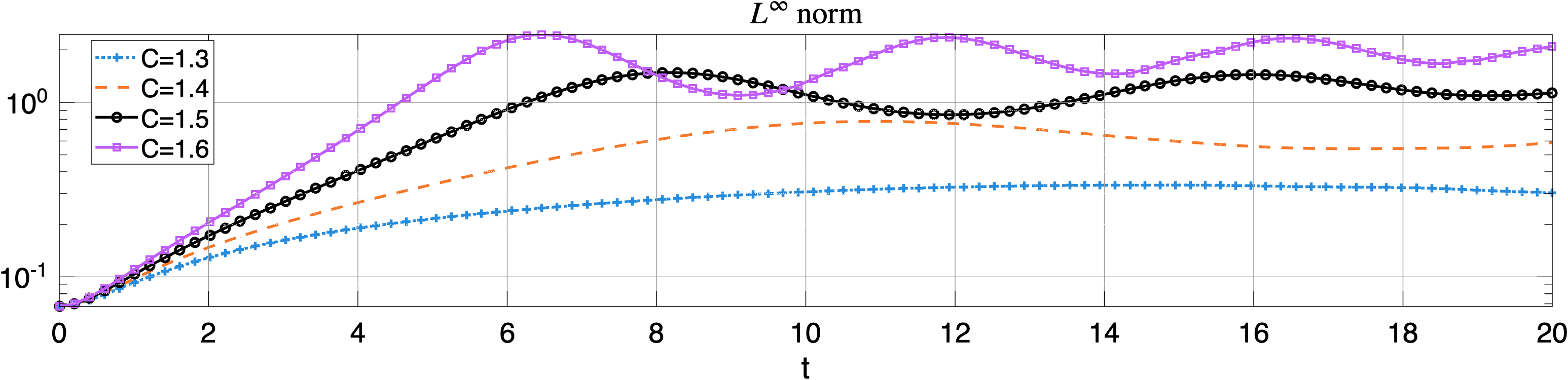}
	\caption{Evolution of the size of $u(t)$ in time for different values of $C.$ The inhomogeneity $u$ satisfies the fully nonlinear equation \eqref{eq:aug1} with all parameters and initial condition as in Figure \ref{fig:3}. {\bf Top:} evolution in time of the $L^2$ norm in space. {\bf Bottom:} evolution in time of the $L^\infty$ norm in space. }
		\label{fig:5}
\end{figure}

In Figure \ref{fig:5} we can see how the qualitative behavior changes with the strength of the background $C.$ The values $C\in[1.3,1.6]$ also lead to instability, but one that grows more slowly and reaches a smaller maximum, compared to the case  $C=1.9.$ 

Clearly, the factor by which $u$ grows is significant,  an indicator both of modulation instability and of the difficulty of the computation. 
One way to measure that, is through the {\em Inhomogeneity Amplification Factor,}
\begin{equation}\label{eq:if}
\mathrm{IAF} :=  \frac{\max\limits_{t,x,y} |u(x,y,t)|}{\max\limits_{x,y} | u_0(x,y)|}.
\end{equation}

In stable problems, the IAF is always very close to $1.$
In problems where modulation instability is present,
 the IAF itself is a function of the initial condition as well as the background $\Gamma.$ Indeed, since the meta-stable pattern has values of $O(1),$ then $u$ will grow as much as needed from its initial value to reach these values, $\mathrm{IAF}\cdot \|u_0\|_{L^\infty} =O(1)$ (numerical experiments verify this). In all the computations here we start with inhomogeneities that are around $5\cdot10^{-2}$  in  $L^\infty$ sense, so that the IAFs are loosely comparable among different examples. 
 
 The main purpose of the  IAF is to capture how challenging the computation is numerically: a computation where the solution grows by a factor of $100$ is more challenging than a computation where the solution grows by a factor of $5.$
If we look at the physical interpretation of the problem, however, the intrinsic measure of how large the solution becomes should be compared to the homogeneous background $\Gamma.$



%
%
%
%
%
%
%

\subsection{Monte Carlo investigation of amplification factors} \label{sec:4.3}
If we  go back to the derivation of equation \eqref{eq:inhomalb1}, the total second-moment $R(x,y,t)= \mathrm{E} \big[ A(x,t) \overline{A(y,t)} \big]
$ of the underlying sea state is decomposed to
\begin{equation}
	R(x,y,t) = u(x,y,t) + \Gamma(x-y).
\end{equation}

The impact of $u$ on the physical problem is through the total second moment $u(x,y,t) + \Gamma(x-y).$  Following \cite{stiassnie2008recurrent,Ribal2013AlberEquation}, we keep track of the {\em Total Amplification Factor}
\begin{equation}\label{eq:af}
\mathrm{TAF} :=  \frac{\max\limits_{t,x,y} |u(x,y,t) + \Gamma(x-y)|}{\max\limits_{x,y} | \Gamma(x-y)|} = \frac{\max\limits_{t,x,y} |u(x,y,t) + \Gamma(x-y)|}{  \Gamma(0)}.
\end{equation}
Note that the quasi-homogeneity assumption means that initially $u$ is negligible compared to $\Gamma;$ thus whether $u_0$ is included or not in the denominator above  makes no meaningful difference.

The TAF controls how large  localized extreme events can be compared to the homogeneous background. This quantity is referred to as ``$\tilde\rho(0,0,\tilde\tau) / \tilde\rho_h(0)$'' in \cite{stiassnie2008recurrent,Ribal2013AlberEquation}, and in \cite{Ribal2013AlberEquation} the rogue wave condition is reported to be equivalent to
\begin{equation}
\mathrm{TAF} > 2.84.
\end{equation}
Thus, in the ocean waves context,  $\mathrm{TAF}=2$ is extreme and $\mathrm{TAF} > 2.84$ is a rogue event.

The TAF is related to the IAF, but the correspondence is not one to one.
By the direct and inverse triangle inequalities we have that
\begin{equation}
\left|1-\frac{\mathrm{IAF} \cdot \|u_0\|_{L^\infty}}{\Gamma(0)}  \right|
\leqslant	\mathrm{TAF} \leqslant 1 + \frac{ \mathrm{IAF} \cdot \|u_0\|_{L^\infty} }{\Gamma(0)}.
\end{equation}
In other words, large IAF is necessary for large TAF -- but in the unstable case, where $\mathrm{IAF}\cdot \|u_0\|_{L^\infty}=O(1)$ it does not follow immediately that large IAF leads to large TAF.


\medskip

\begin{figure}
	\includegraphics[width=0.9\textwidth]{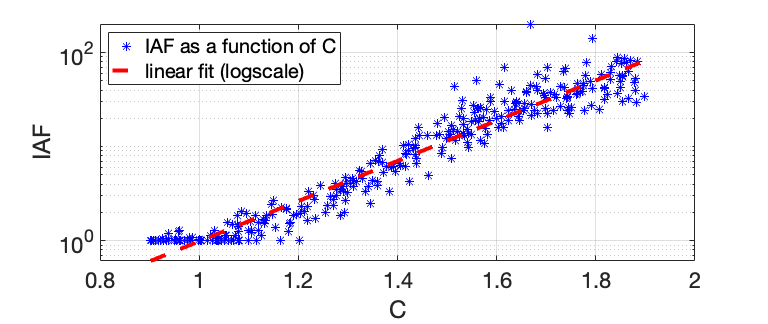}
	\caption{IAF versus intensity of the background, $C$ for randomized initial conditions and different values of $C.$ Details of the computation are given in \eqref{eq:randomized} onwards. The IAF is defined in \eqref{eq:if}. The onset of instability around $C=1$ is clear. Final time is $T=16.$}
	\label{fig:AFs1}
\end{figure}

\begin{figure}
	\includegraphics[width=0.9\textwidth]{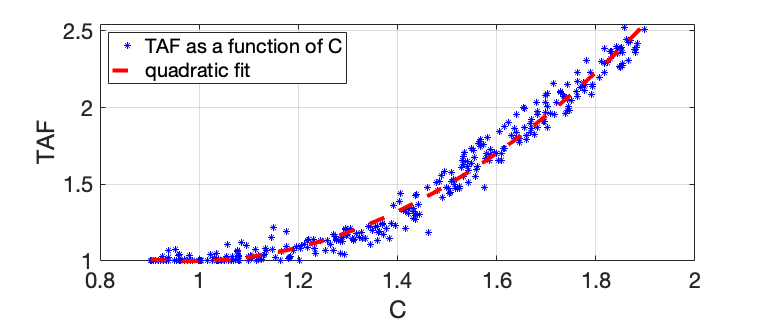}
	\caption{TAF versus intensity of the background, $C$ for randomized initial conditions and different values of $C.$ Details of the computation are given in \eqref{eq:randomized} onwards. The TAF is defined in \eqref{eq:af}. The quadratic dependence on $C$ means that the TAF exhibits no meaningful change at $C=1,$ with the change being noticeable at $C\approx 1.35$ onwards. Final time is $T=16.$}
	\label{fig:AFs2}
\end{figure}

\begin{figure}
	\includegraphics[width=0.9\textwidth]{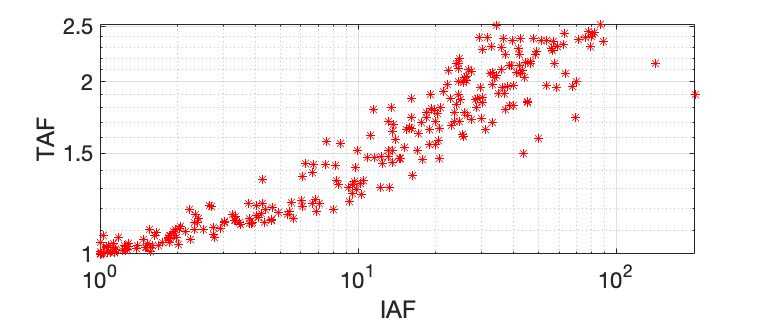}
	\caption{Scatter diagram of TAF versus IAF.}
	\label{fig:AFs3}
\end{figure}

\begin{figure}
	\includegraphics[width=0.9\textwidth]{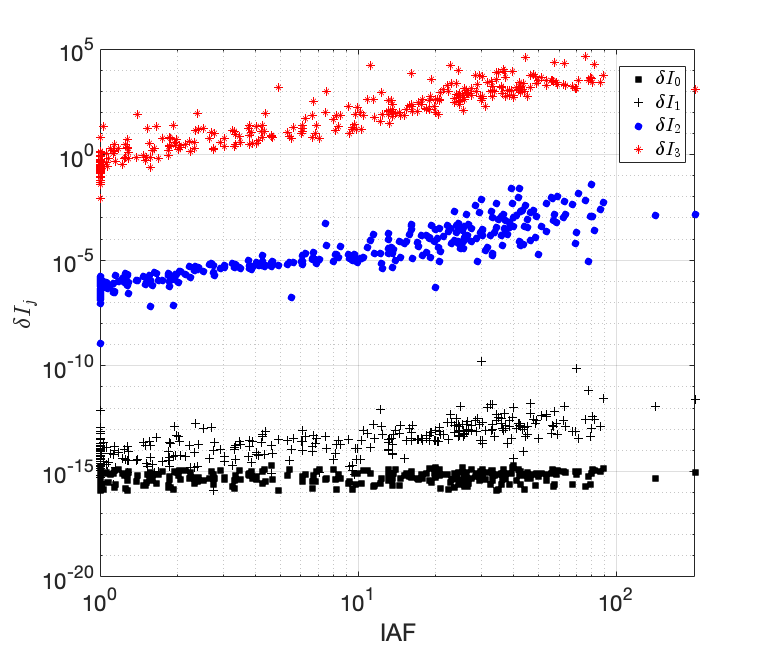}
	\caption{Relative errors in the observables $I_j$ across all simulations, plotted against IAF.  $I_0,I_1$ are conserved to 12 digits in all cases. $I_2$ is approximated to between 5 and 2 digits, with larger IAF leading to worse approximation. $I_3,$ which includes second derivatives, can only be approximated to 2 digits for stable problems, i.e., when $\mathrm{IAF}=1.$ Once the solution grows substantially, $I_3$ can no longer be approximated. }
	\label{fig:AFs4}
\end{figure}

To see how the IAF and the TAF behave for various intensities of the background $C,$ and various randomized initial conditions (all of comparable $L^2$ and $L^\infty$ sizes), we perform a Monte Carlo investigation. In what follows we compute amplification factors IAF and TAF for different values of $C$ in \eqref{eq:gamma}, and different values of the coefficients $A_j$  in equation  \eqref{eq:u0}. More specifically, we set
\begin{equation}\label{eq:randomized}
C \sim \mathcal{U}([0.9,1.9]), \qquad 
	A_j \sim \mathcal{N}(0,\frac{1}9) + i \mathcal{N}(0,\frac{1}9).
\end{equation}
The resulting problem of the form \eqref{eq:aug1} is then solved for $t\in [0,16],$ on $x,y\in[-L/2,L/2]^2$ with $L=50.$ This is long enough time for the maximum values to be attained, according to testing. Discretization parameters are $\tau = 2\cdot 10^{-3},$ $h=0.12.$ This is somewhat coarser resolution than in Section \ref{subsec:42} for reasons of computational performance. The resulting amplification factors, along with quality control metadata, are plotted in Figures \ref{fig:AFs1}-\ref{fig:AFs4}.

The numerical simulations for the Monte Carlo investigation were performed on the standard nodes of Hypatia HPC cluster at the University of St Andrews (each standard node contains two  84-core AMD EPYC 9634 CPUs, 384GB of memory and a local 1.7TB scratch disk). We used MATLAB R2025b with SLURM for job scheduling. We executed $310$ independent realizations of the PDE initial value problem solver in parallel batches, each allocated 6 CPU cores to maximize BLAS efficiency.

Figure \ref{fig:AFs1} clearly indicates that $u$ starts to grow noticeably from its initial size around $C\approx 1,$ in agreement with the linear stability analysis. Figure \ref{fig:AFs2} shows that the maximum value of the coherent structure depends clearly on the background intensity $C.$ The quadratic nature of this dependence means that for $C$ close to 1 the TAF is not noticeably larger than $1.$ On the other hand, for $C\approx 1.4,$ we have $\mathrm{TAF}\approx 1.4.$ This means that the variance of the envelope amplitude is $40\%$ higher than what it would be in the purely homogeneous case, a substantial increase.

Figure \ref{fig:AFs3} shows that the TAF and IAF are positively correlated, but still there is substantial room in their values, i.e., they are not in a one-to-one relationship. Figure \ref{fig:AFs4} is the quality control for this sequence of runs: even with solution growing hundreds of times from their initial size, and relatively long final time of $T=16,$ $I_0,I_1$ are preserved to 12 digits or more, and $I_2$ is approximated to at least 2 digits. On the other hand, $I_3$ cannot be well approximated for long times when the solution grows substantially.

It is clearly seen that the amplification factors are not sensitive to the exact details of the initial inhomogeneity; moreover, a robust dependence on $C$ clearly emerges. To the best of our knowledge, this clear dependence of the maxima of the fully developed modulation instability on the background strength, and lack of substantial dependence on the initial condition, was not known before.

\section{Conclusions} \label{sec:Conclusions}

\subsection{The  numerical method}
The von Neumann equation with singular interaction is widely used as a model for the second moments of stochastic ocean wave fields. As discussed above, only one numerical method existed for this equation prior to the present work \cite{stiassnie2008recurrent,Ribal2013AlberEquation}.

We introduce
a relaxation-Crank-Nicolson scheme which is second-order in time and requires the solution of a linear system in each timestep, i.e., it is explicit in the nonlinearity. This is achieved by the introduction of an auxiliary variable on a staggered time-grid. 
The scheme is implemented with fourth-order finite differences in space. It inherits at the  discrete level the problem's $L^2$ balance law (Proposition \ref{prop:discrconbll}) and two invariants, making it  particularly well suited to investigate the possible growth of the solution. We also show consistency and boundedness results for the scheme, while noting that full global stability still requires stronger a priori results than are available for this problem. 

In Section \ref{sec:2} the scheme is validated with a soliton-derived exact solution. It is found that the initialization of the auxiliary variable is important for the scheme to  achieve the expected order in both space and time, for both the main and auxiliary variables. Moreover, it is found that the main variable is to a surprising extent not affected by the error in the auxiliary variable, for example exhibiting the correct order of convergence even if the auxiliary variable does not (Tables \ref{eoch6} and \ref{eoch2}). This seems to indicate an averaging effect, and fully understanding it could be vital  for the analysis of this type of scheme. Finally,  Figure \ref{Fig:4} shows that the bulk of the error is  due to error in the phase velocity. This is to be expected from the von Neumann stability analysis for the linear part of the equation (where our method becomes the Pad\'e(2,2) method).

\subsection{New findings on the bifurcation}
The dispersion dominated regime (cf. Figure \ref{fig:8}) is clearly recovered, and it is controlled by the linear theory. In the unstable regime, the linear stability analysis accurately predicts the initial growth phase, but fails to forecast the maximum amplitude. In the fully developed, nonlinear phase of the instability, a meta-stable pattern is found, cf. Figure \ref{fig:3}. This includes the recurrent maxima reported in \cite{Ribal2013AlberEquation}, but it is not just one recurrent maximum: the use of the larger spatial domain allows us to document for the first time the full space-time lattice. Its main scalings (e.g., maximum value and spacing of maxima) seem to depend mainly on the homogeneous background and not on the initial condition, according to extensive numerical investigation. This meta-stable pattern is not well understood to date, e.g., there are no explicit or asymptotic approximations for it.

To quantify the strength of the instability and its physical consequences, we use non-dimensional amplification factors.
 The TAF, inspired by \cite{Ribal2013AlberEquation} and defined in \eqref{eq:af}, measures the increase in the second moment of the whole sea state due to the modulation instability. Physically, this controls the variance of the envelope amplitude -- thus a TAF of $1.4$ means $40\%$ higher variance at the local maxima of $u.$ On the other hand, the IAF, introduced in \eqref{eq:if}, measures how much the inhomogeneity grows. A large value of IAF, say $\mathrm{IAF}=7,$  is challenging numerically, and also indicative of the presence of modulation instability. On the other hand, large values of IAF are not sufficient by themselves to produce noticeably large TAF (say, $\mathrm{TAF}>1.1$).
It is found for the first time that, for a given spectral shape, both the IAF and TAF have a clear dependence on the intensity of the background and very little dependence on the initial inhomogeneity.

The study of the two kinds of amplification factors indicates that onset of modulation instability (i.e., large IAF) does not automatically guarantee large extreme events (i.e., large TAF) --  much less rogue waves. On the other hand, sufficiently strong instability does guarantee substantially larger variance of the sea surface, and thus larger extreme events.

If the variance becomes very large  (e.g., the $\mathrm{TAF}>2.84$ criterion of \cite{Ribal2013AlberEquation}) then it effectively guarantees a rogue wave.  However, this would need a very strong instability, almost twice the intensity required for onset of modulation instability. Such sea states seem to exist on the record, but just barely \cite{Athanassoulis2018}. Moreover, the recurrent character of the maxima (cf. Figure \ref{fig:3}) would point towards a persistent series of rogue waves, not a single event.

Thus, it seems more likely  that most rogue waves appear for $\mathrm{TAF}<2.84.$ In a modulationally unstable sea there are recurrent hotspots of increased variance -- cf. the peaks of Figure \ref{fig:3}. Intensities $30-40\%$ higher than what is required for the onset of instability  could give rise to $30-40\%$ higher variance over the hotspots, according to Figure \ref{fig:AFs2}.  These recurrent hotspots now don't necessarily guarantee rogue waves, but they do make them much more likely. This suggests a possible mechanism by which rogue waves become orders of magnitude more likely than predicted by linear statistics in intense, modulationally unstable sea states. At the same time, such events remain rare and do not exhibit persistent recurrence.

%
%

\subsection{Code availability}

A standalone version of the solver, including implementations of many of the features discussed here (EOC, constraint error, invariants, growth of inhomogeneity, Monte Carlo investigation of the amplification factor) can be found in \url{https://github.com/aathanas/rcn_Alber_MATLAB}.

\bibliographystyle{siam}
\bibliography{CollectionSupervised.bib}

\begin{thebibliography}{10}

\bibitem{Alber1978}
{\sc I.~Alber}, {\em The effects of randomness on the stability of two-dimensional surface wavetrains}, Proceedings of the Royal Society of London. A. Mathematical and Physical Sciences, 363 (1978), pp.~525--546.

\bibitem{Alsafri2023}
{\sc N.~Alsafri}, {\em {Numerical Methods For A Model Of Two-Phase Flows}}, {PhD thesis}, University of Dundee, 2023.

\bibitem{Arnold1994AbsorbingBC}
{\sc A.~Arnold}, {\em On absorbing boundary conditions for quantum transport equations}, ESAIM: Mathematical Modelling and Numerical Analysis, 28 (1994), pp.~853--872.

\bibitem{arnold1992numerical}
{\sc A.~Arnold and F.~Nier}, {\em {Numerical analysis of the deterministic particle method applied to the Wigner equation}}, Mathematics of Computation, 58 (1992), pp.~645--669.

\bibitem{Arnold1996WignerPoisson}
{\sc A.~Arnold and C.~Ringhofer}, {\em An operator splitting method for the wigner--poisson problem}, SIAM Journal on Numerical Analysis, 33 (1996), pp.~1622--1643.

\bibitem{Athanassoulis2018}
{\sc A.~Athanassoulis, G.~Athanassoulis, M.~Ptashnyk, and T.~Sapsis}, {\em {Strong solutions for the Alber equation and stability of unidirectional wave spectra}}, Kinetic and Related Models, 13 (2020), pp.~703--737.

\bibitem{athanassoulis2024efficient}
{\sc A.~Athanassoulis, T.~Katsaounis, and I.~Kyza}, {\em Efficient numerical approximations for a nonconservative nonlinear schr{\"o}dinger equation appearing in wind-forced ocean waves}, Studies in Applied Mathematics, 153 (2024), p.~e12774.

\bibitem{Athanassoulis2023a}
{\sc A.~Athanassoulis, T.~Katsaounis, I.~Kyza, and S.~Metcalfe}, {\em {A novel, structure-preserving, second-order-in-time relaxation scheme for Schr{\"{o}}dinger-Poisson systems}}, Journal of Computational Physics, 490 (2023), p.~112307.

\bibitem{athanassoulis2011strong}
{\sc A.~Athanassoulis and T.~Paul}, {\em Strong phase-space semiclassical asymptotics}, SIAM journal on mathematical analysis, 43 (2011), pp.~2116--2149.

\bibitem{Athanassoulis2017RogueWaves}
{\sc A.~G. Athanassoulis, G.~A. Athanassoulis, and T.~Sapsis}, {\em {Localized instabilities of the Wigner equation as a model for the emergence of Rogue Waves}}, Journal of Ocean Engineering and Marine Energy, 3 (2017), pp.~325--338.

\bibitem{Athanassoulis2023b}
{\sc A.~G. Athanassoulis and I.~Kyza}, {\em {Modulation instability and convergence of the random phase approximation for stochastic sea states}}, Water Waves, 6 (2024), pp.~145--167.

\bibitem{benjamin1967disintegration}
{\sc T.~B. Benjamin and J.~E. Feir}, {\em The disintegration of wave trains on deep water part 1. theory}, Journal of Fluid Mechanics, 27 (1967), pp.~417--430.

\bibitem{besse2004relaxation}
{\sc C.~Besse}, {\em {A relaxation scheme for the nonlinear Schr{\"o}dinger equation}}, SIAM Journal on Numerical Analysis, 42 (2004), pp.~934--952.

\bibitem{besse2021energy}
{\sc C.~Besse, S.~Descombes, G.~Dujardin, and I.~Lacroix-Violet}, {\em {Energy-preserving methods for nonlinear Schr{\"o}dinger equations}}, IMA Journal of Numerical Analysis, 41 (2021), pp.~618--653.

\bibitem{biondini2018universal}
{\sc G.~Biondini, S.~Li, D.~Mantzavinos, and S.~Trillo}, {\em Universal behavior of modulationally unstable media}, SIAM Review, 60 (2018), pp.~888--908.

\bibitem{chaub2025semiclassical}
{\sc T.~Chaub}, {\em Semiclassical limit and singular Vlasov Equations}, PhD thesis, Universit{\'e} Paris-Saclay, 2025.

\bibitem{chen2017global}
{\sc T.~Chen, Y.~Hong, and N.~Pavlovi{\'c}}, {\em Global well-posedness of the nls system for infinitely many fermions}, Archive for rational mechanics and analysis, 224 (2017), pp.~91--123.

\bibitem{Chen2022SpectralWignerPoisson}
{\sc Z.~Chen, H.~Jiang, and S.~Shao}, {\em A higher-order accurate operator splitting spectral method for the wigner--poisson system}, Journal of Computational Electronics, 21 (2022), pp.~1313--1329.

\bibitem{Chen2019Wigner4D}
{\sc Z.~Chen, S.~Shao, and W.~Cai}, {\em {A high order efficient numerical method for 4-D Wigner equation of quantum double-slit interferences}}, Journal of Computational Physics, 396 (2019), pp.~54--71.

\bibitem{Chen2019UnboundedWigner}
{\sc Z.~Chen, Y.~Xiong, and S.~Shao}, {\em Numerical methods for the wigner equation with unbounded potential}, Journal of Scientific Computing, 79 (2019), pp.~1275--1305.

\bibitem{chong2024many}
{\sc J.~J. Chong, L.~Lafleche, and C.~Saffirio}, {\em From many-body quantum dynamics to the hartree-fock and vlasov equations with singular potentials.}, Journal of the European Mathematical Society (EMS Publishing), 26 (2024).

\bibitem{dematteis2018rogue}
{\sc G.~Dematteis, T.~Grafke, and E.~Vanden-Eijnden}, {\em Rogue waves and large deviations in deep sea}, Proceedings of the National Academy of Sciences, 115 (2018), pp.~855--860.

\bibitem{Dorda2015WENO}
{\sc A.~Dorda and F.~Sch{\"u}rrer}, {\em {A WENO-solver combined with adaptive momentum discretization for the Wigner transport equation and its application to resonant tunneling diodes}}, Journal of Computational Physics, 284 (2015), pp.~95--113.

\bibitem{dysthe2008oceanic}
{\sc K.~Dysthe, H.~E. Krogstad, and P.~M{\"u}ller}, {\em Oceanic rogue waves}, Annual Review of Fluid Mechanics, 40 (2008), pp.~287--310.

\bibitem{Filbet2025SemiClassicalVN}
{\sc F.~Filbet and F.~Golse}, {\em On the approximation of the von-neumann equation in the semi-classical limit. part i: Numerical algorithm}, Journal of Computational Physics, 488 (2025), p.~112336.

\bibitem{Gramstad2017}
{\sc O.~Gramstad}, {\em {Modulational Instability in JONSWAP Sea States Using the Alber Equation}}, in ASME 2017 36th International Conference on Ocean, Offshore and Arctic Engineering, 2017.

\bibitem{grekhneva2020dynamics}
{\sc A.~D. Grekhneva and V.~Z. Sakbaev}, {\em {Dynamics of a set of quantum states generated by a nonlinear Liouville--von Neumann equation}}, Computational Mathematics and Mathematical Physics, 60 (2020), pp.~1337--1347.

\bibitem{hadama2025global}
{\sc S.~Hadama and Y.~Hong}, {\em Global well-posedness of the nonlinear hartree equation for infinitely many particles with singular interaction}, Journal of Functional Analysis,  (2025), p.~111102.

\bibitem{han2025semiclassical}
{\sc D.~Han-Kwan and F.~Rousset}, {\em Semiclassical limit of cubic nonlinear schr\"odinger equations for mixed states}, arXiv preprint arXiv:2510.21313,  (2025).

\bibitem{Jiang2023WignerPoissonRTD}
{\sc H.~Jiang, T.~Lu, W.~Yao, and W.~Zhang}, {\em Numerical study of transient wigner--poisson model for rtds: Numerical method and its applications}, SIAM Journal on Scientific Computing, 45 (2023), pp.~B1203--B1228.

\bibitem{Jiang2021HybridWigner}
{\sc H.~Jiang, T.~Lu, and X.~Yin}, {\em A hybrid explicit-implicit scheme for the time-dependent wigner equation}, Journal of Computational Mathematics, 39 (2021), pp.~22--44.

\bibitem{lions1993mesures}
{\sc P.-L. Lions and T.~Paul}, {\em Sur les mesures de wigner}, Revista matem{\'a}tica iberoamericana, 9 (1993), pp.~553--618.

\bibitem{mei2005theory}
{\sc C.~C. Mei, M.~A. Stiassnie, and D.~K.-P. Yue}, {\em {Theory and applications of ocean surface waves: Part 2: nonlinear aspects}}, World Scientific, 2005.

\bibitem{Muscato2016StochasticWigner}
{\sc O.~Muscato and W.~Wagner}, {\em A class of stochastic algorithms for the wigner equation}, SIAM Journal on Scientific Computing, 38 (2016), pp.~B772--B796.

\bibitem{Muscato2019StochasticWigner}
\leavevmode\vrule height 2pt depth -1.6pt width 23pt, {\em A stochastic algorithm without time discretization error for the wigner equation}, Kinetic and Related Models, 12 (2019), pp.~1--23.

\bibitem{neumann1927wahrscheinlichkeitstheoretischer}
{\sc J.~v. Neumann}, {\em Wahrscheinlichkeitstheoretischer aufbau der quantenmechanik}, Nachrichten von der Gesellschaft der Wissenschaften zu G{\"o}ttingen, Mathematisch-Physikalische Klasse, 1927 (1927), pp.~245--272.

\bibitem{ochi1998ocean}
{\sc M.~K. Ochi}, {\em Ocean waves}, Cambridge University Press, 1998.

\bibitem{Onorato2003LandauDamping}
{\sc M.~Onorato, A.~Osborne, R.~Fedele, and M.~Serio}, {\em Landau damping and coherent structures in narrow-banded deep water gravity waves}, Physical Review E, 67 (2003), p.~046305.

\bibitem{onorato2009statistical}
{\sc M.~Onorato, T.~Waseda, A.~Toffoli, L.~Cavaleri, O.~Gramstad, P.~Janssen, T.~Kinoshita, J.~Monbaliu, N.~Mori, A.~R. Osborne, et~al.}, {\em {Statistical Properties of Directional Ocean Waves: The Role of the Modulational Instability in the Formation of Extreme Events}}, Physical Review Letters, 102 (2009), p.~114502.

\bibitem{Ribal2013AlberEquation}
{\sc A.~Ribal, A.~V. Babanin, I.~Young, A.~Toffoli, and M.~Stiassnie}, {\em {Recurrent solutions of the Alber equation initialized by Joint North Sea Wave Project spectra}}, Journal of Fluid Mechanics, 719 (2013), pp.~314--344.

\bibitem{ringhofer1992spectral}
{\sc C.~Ringhofer}, {\em A spectral collocation technique for the solution of the wigner--poisson problem}, SIAM journal on numerical analysis, 29 (1992), pp.~679--700.

\bibitem{Ringhofer1989AbsorbingBC}
{\sc C.~Ringhofer, D.~Ferry, and N.~Kluksdahl}, {\em Absorbing boundary conditions for the simulation of quantum transport phenomena}, Transport Theory and Statistical Physics, 18 (1989), pp.~605--618.

\bibitem{rosati2013wigner}
{\sc R.~Rosati, F.~Dolcini, R.~C. Iotti, and F.~Rossi}, {\em Wigner-function formalism applied to semiconductor quantum devices: Failure of the conventional boundary condition scheme}, Physical Review B---Condensed Matter and Materials Physics, 88 (2013), p.~035401.

\bibitem{sakbaev2022blow}
{\sc V.~Z. Sakbaev and A.~Shiryaeva}, {\em Blow-up of states in the dynamics given by the schr{\"o}dinger equation with a power-law nonlinearity in the potential}, Differential Equations, 58 (2022), pp.~497--508.

\bibitem{sapsis2021statistics}
{\sc T.~P. Sapsis}, {\em Statistics of extreme events in fluid flows and waves}, Annual Review of Fluid Mechanics, 53 (2021), pp.~85--111.

\bibitem{Schulz2021ExponentialIntegrators}
{\sc L.~Schulz, B.~Inci, M.~Pech, and D.~Schulz}, {\em Subdomain-based exponential integrators for quantum liouville-type equations}, Journal of Computational Electronics, 20 (2021), pp.~2027--2038.

\bibitem{Shukla2007ModulationalIncoherent}
{\sc P.~K. Shukla, M.~Marklund, and L.~Stenflo}, {\em Modulational instability of nonlinearly interacting incoherent sea states}, JETP Letters, 85 (2007), pp.~33--36.

\bibitem{stiassnie2008recurrent}
{\sc M.~Stiassnie, A.~Regev, and Y.~Agnon}, {\em {Recurrent solutions of Alber's equation for random water-wave fields}}, Journal of Fluid Mechanics, 598 (2008), pp.~245--266.

\bibitem{Tian2012LvNFiniteTemp}
{\sc H.~Tian and G.~H. Chen}, {\em An efficient solution of liouville-von neumann equation that is applicable to zero and finite temperatures}, The Journal of Chemical Physics, 137 (2012), p.~204114.

\bibitem{weideman1995computing}
{\sc J.~Weideman}, {\em Computing the {Hilbert} transform on the real line}, Mathematics of Computation, 64 (1995), pp.~745--762.

\bibitem{wigner1932quantum}
{\sc E.~Wigner}, {\em On the quantum correction for thermodynamic equilibrium}, Physical review, 40 (1932), p.~749.

\bibitem{zakharov1968stability}
{\sc V.~E. Zakharov}, {\em Stability of periodic waves of finite amplitude on the surface of a deep fluid}, Journal of Applied Mechanics and Technical Physics, 9 (1968), pp.~190--194.

\bibitem{zouraris2021error}
{\sc G.~E. Zouraris}, {\em Error estimation of the besse relaxation scheme for a semilinear heat equation}, ESAIM: Mathematical Modelling and Numerical Analysis, 55 (2021), pp.~301--328.

\bibitem{zouraris2023error}
\leavevmode\vrule height 2pt depth -1.6pt width 23pt, {\em {Error estimation of the relaxation finite difference scheme for the nonlinear Schr{\"o}dinger equation}}, SIAM Journal on Numerical Analysis, 61 (2023), pp.~365--397.

\end{thebibliography}

\appendix

\section{Local well-posedness for the continuous problem}\label{App:A}

\noindent {\bf Proof of Theorem \ref{thrm:torusreg}:}
To control the trace nonlinearity, we will use the Wiener algebra norm, defined as
\begin{equation}
	\tnorm{u} := \sum\limits_{k,l} |\hat{u}_{k,l}|, \qquad \tnorm{\Gamma} := \sum\limits_{n\in\mathbb{Z}} |P_n|.
\end{equation}
Obviously, $\|u\|_{L\infty} \leqslant \tnorm{u};$ moreover, 
the Wiener algebra norm allows the bound
\begin{equation}
	\tnorm{ \phi(t) u(t) } = 	\tnorm{ (u(x,x,t) - u(y,y,t)) u(x,y,t) } \leqslant 	2\tnorm{u(t)}^2.
\end{equation}
This will be a crucial ingredient for the existence of solutions.

Let us 
 rewrite equation \eqref{eq:aug1} using the Fourier series:
\[
\begin{array}{c}
i \sum\limits_{k,l \in \mathbb{Z}} \partial_t \hat{u}_{k,l}(t) e^{2\pi i \frac{kx+ly}L} + p \sum\limits_{k,l \in \mathbb{Z}} \hat{u}_{k,l}(t) \left((2\pi i \frac{k}L)^2 - (2\pi i \frac{l}L)^2\right) e^{2\pi i \frac{kx+ly}L}+ \\
+q  \left(\sum\limits_{n \in \mathbb{Z}} P_n e^{2\pi i \frac{n(x-y)}L} + \sum\limits_{k,l \in \mathbb{Z}} \hat{u}_{k,l}(t) e^{2\pi i \frac{kx+ly}L}  \right)
\left(  \sum\limits_{k,l \in \mathbb{Z}} \hat{u}_{k,l}(t) \left[e^{2\pi i \frac{k+l}L x} - e^{2\pi i \frac{k+l}L y} \right]\right)  =0. 
\end{array}
\]

Take the inner product with $e^{2\pi i \frac{k'x+l'y}L};$ for clarity we work out explicitly the last two terms term:
\[
\begin{array}{c}
	q\bigintsss \left(\sum\limits_{n \in \mathbb{Z}} P_n e^{2\pi i \frac{n(x-y)}L}\right)
\left(  \sum\limits_{k,l \in \mathbb{Z}} \hat{u}_{k,l}(t) \left[e^{2\pi i \frac{k+l}L x} - e^{2\pi i \frac{k+l}L y} \right]\right) (e^{-2\pi i \frac{k'x+l'y}L}) dxdy=\\

=q\sum\limits_{n,k,l}P_n \hat{u}_{k,l}(t) \bigintsss \left\{    e^{2\pi i [x\frac{n+k+l-k'}L + y \frac{-n-l'}L]}
- e^{2\pi i [x\frac{n-k'}L + y \frac{-n+k+l-l'}L]}\right\} dxdy=\\

=qL^2\sum\limits_{n,k,l}P_n \hat{u}_{k,l}(t)   \left[  \delta_{n+k+l,k'} \delta_{n+l'} - \delta_{n,k'}\delta_{k+l-n,l'}  \right]=\\

=qL^2\left[ \sum\limits_{k,l}P_{-l'} \hat{u}_{k,l}(t) \delta_{k+l,k'+l'}  -\sum\limits_{k,l}P_{k'} \hat{u}_{k,l}(t)  \delta_{k+l,k'+l'}  \right]=\\
=qL^2[P_{-l'}-P_{k'}] \sum\limits_{k+l=k'+l'} \hat{u}_{k,l}(t) = qL^2[P_{-l'}-P_{k'}] \sum\limits_{S \in \mathbb{Z}} \hat{u}_{S,k'+l'-S}(t) ,
\end{array}
\]
and
\[
\begin{array}{c}
	q\bigintsss \left( \sum\limits_{a,b \in \mathbb{Z}} \hat{u}_{a,b}(t) e^{2\pi i \frac{ax+by}L}  \right)
\left(  \sum\limits_{k,l \in \mathbb{Z}} \hat{u}_{k,l}(t) \left[e^{2\pi i \frac{k+l}L x} - e^{2\pi i \frac{k+l}L y} \right]\right) (e^{-2\pi i \frac{k'x+l'y}L}) dxdy=\\

=q\sum\limits_{a,b,k,l} \hat{u}_{a,b}(t) \hat{u}_{k,l}(t) \bigintsss \left\{    e^{2\pi i [x\frac{a+k+l-k'}L + y \frac{b-l'}L]}
- e^{2\pi i [x\frac{a-k'}L + y \frac{b+k+l-l'}L]}\right\} dxdy=\\

=qL^2\sum\limits_{a,b,k,l} \hat{u}_{a,b}(t) \hat{u}_{k,l}(t)  \left[  \delta_{a+k+l,k'} \delta_{b,l'} - \delta_{a,k'}\delta_{b+k+l,l'}  \right]=\\

=qL^2\left( \sum\limits_{a,k,l} \hat{u}_{a,l'}(t) \hat{u}_{k,l}(t)    \delta_{a+k+l,k'}  - \sum\limits_{b,k,l} \hat{u}_{k',b}(t) \hat{u}_{k,l}(t) \delta_{b+k+l,l'}  \right)=\\

=qL^2\left( \sum\limits_{K,\Lambda} \hat{u}_{k'-K-\Lambda,l'}(t) \hat{u}_{K,\Lambda}(t)      - \sum\limits_{K,\Lambda} \hat{u}_{k',l'-K-\Lambda}(t) \hat{u}_{K,\Lambda}(t)  \right).
\end{array}
\]

Finally if we multiply by $-i/L^2$ and suppress the primes we obtain
\begin{equation}\label{eq:poikjuyh}
\begin{array}{c}
  \partial_t\hat{u}_{k,l}(t)  +i \frac{4\pi^2 p}{L^2}  \hat{u}_{k,l}(t) (k-l)(k+l) 
-iq \left[  P_{-l} -  P_{k} \right]  \sum\limits_{S\in \mathbb{Z}}   \hat{u}_{S,k+l-S}(t) \\
-iq  \left( \sum\limits_{K,\Lambda} \hat{u}_{k-K-\Lambda,l}(t) \hat{u}_{K,\Lambda}(t)      - \sum\limits_{K,\Lambda} \hat{u}_{k,l-K-\Lambda}(t) \hat{u}_{K,\Lambda}(t)  \right)
=0.
\end{array}
\end{equation}

Integrating in time leads to
\[
\begin{array}{c}
\hat{u}_{k,l}(t) = e^{-i \frac{4\pi^2 p}{L^2}  (k-l)(k+l) t} \hat{u}_{k,l}(0) +   
iq \int\limits_{\tau=0}^t  e^{-i \frac{4\pi^2 p}{L^2}  (k^2-l^2) (t-\tau)} \left[  P_{-l} -  P_{k} \right]  \sum\limits_{S\in \mathbb{Z}}   \hat{u}_{S,k+l-S}(\tau) d\tau +\\
+ iq \int\limits_{\tau=0}^t  e^{-i \frac{4\pi^2 p}{L^2}  (k^2-l^2) (t-\tau)} 
\sum\limits_{K,\Lambda}\Bigl( \left(  \hat{u}_{k-K-\Lambda,l}(\tau)       -  \hat{u}_{k,l-K-\Lambda}(\tau)   \right) \hat{u}_{K,\Lambda}(\tau) \Bigr)
d\tau.
\end{array}
\]
By taking the $l^1$ norm in $k,l$ it follows that
\begin{equation}
	\tnorm{u(t)} \leqslant 	\tnorm{u(0)} + 2 |q|  \int\limits_{\tau=0}^t \left( \tnorm{\Gamma} \cdot \tnorm{u(\tau)} + \tnorm{u(\tau)}^2 \right) d\tau.
\end{equation}
This is a standard bootstrap argument leading to local-in-time existence and uniqueness; one easily checks that $\tnorm{u(t)}<\infty$  at least for $0 \leqslant t < T_*$ where
\[
T_*=
\begin{cases}
\dfrac{1}{2 |q| \tnorm{\Gamma}}
\ln \left( 1+ \dfrac{\tnorm{\Gamma}}{ \tnorm{u(0)} }\right),
& \text{if } \tnorm{\Gamma}>0, \\[16pt]
\dfrac{1}{2 q \tnorm{u(0)}},
& \text{if }  \tnorm{\Gamma}=0 .
\end{cases}
\]
In particular, denoting by $\Sigma$ the Wiener algebra on the torus, it follows that  problem \ref{eq:aug1} has a continuous propagator
\begin{equation}
U(t): \Sigma \to \Sigma \qquad \mbox{ for } t \in [0,T_*).
\end{equation} 

Now if we apply $\partial_x^a\partial_y^b$ to problem \eqref{eq:aug1}, we get the system
\begin{equation}\label{eq:augk}
\begin{array}{c}
\displaystyle  \phi(x,y,t) = u(x,x,t) - u(y,y,t), \\[4pt]
\displaystyle  i \partial_t (\partial_x^a\partial_y^b u) + p (\Delta_x - \Delta_y)  (\partial_x^a\partial_y^b u) + q \phi(x,y,t) \,  \Bigl(\partial_x^a\partial_y^b\Gamma(x-y) + \partial_x^a\partial_y^b u \Bigr)  =  \\
 \qquad \qquad \qquad \qquad \qquad   = -q\sum\limits_{\nu+\mu < a+b} \binom{a}{\nu} \binom{b}{\mu} \Bigl(\partial_x^{a-\nu}\partial_y^{b-\mu} \phi \Bigr) \Bigl(\partial_x^\nu\partial_y^\mu\Gamma(x-y) + \partial_x^\nu\partial_y^\mu u \Bigr)
\end{array}
\end{equation}
(using standard multi-index notation). Denoting
\begin{equation}
	\tnorm{u}_s = \max\limits_{|a+b|\leqslant s} \tnorm{\partial_x^a\partial_y^b u},
\end{equation}
the mild form of equation \eqref{eq:augk} leads to
\begin{equation}
	\tnorm{u(t)}_s \leqslant \tnorm{u(0)}_s + C \int\limits_{\tau =0}^t 	\tnorm{u(\tau)}_{s-1} 	\tnorm{u(\tau)}_s d\tau.
\end{equation}
Thus, $s$-degree regularity in the Wiener algebra sense can be propagated up to time $T_*$ by working recursively in $s.$

Moreover, if we apply $\partial_t^m$ to problem \eqref{eq:aug1}, we get 
\begin{equation}\label{eq:augm}
\begin{array}{c}
\displaystyle  \phi(x,y,t) = u(x,x,t) - u(y,y,t), \\[4pt]
\displaystyle  i \partial_t (\partial_t^m u) + p (\Delta_x - \Delta_y)  (\partial_t^m u) + q \phi(x,y,t) \,  \partial_t^m u  =  \qquad \qquad  \\
 \qquad \qquad \qquad \qquad \qquad   = -q\sum\limits_{\mu < m} \binom{m}{\mu} \Bigl(\partial_t^{m-\mu}\phi \Bigr) \Bigl(\partial_t^\mu \Gamma(x-y) + \partial_t^\mu u \Bigr)
\end{array}
\end{equation}
(where of course $\partial_t^\mu \Gamma(x-y) \neq 0$ only for $\mu=0$). Now the left hand-side is an Alber equation with $\Gamma=0.$ It is still covered by the previous theory as a special case, and it has a propagator in $\Sigma$ which, by an abuse of notation we will still denote by $U(t).$ Since $\tnorm{ (\partial_t^\mu \phi ) \, (\partial_t^m u) } \leqslant \tnorm{\partial_t^\mu u} \tnorm{ \partial_t^m u},$ the mild form of \eqref{eq:augm} now leads to
\begin{equation}
		\tnorm{\partial_t^m u(t)} \leqslant \tnorm{\partial_t^m u(0)}_s + C \int\limits_{\tau =0}^t 	\tnorm{\partial_t^m u(\tau)} 	\max\limits_{\mu < m}\tnorm{\partial_t^\mu u(\tau)} d\tau.
\end{equation}
By setting $t=0$ and solving for $\partial_t u$ in \eqref{eq:aug1} it follows that
\[
\tnorm{\partial_t u(0)} \leqslant C \tnorm{ u(0)}_2.
\]
More generally, for higher order derivatives in time, one forms a system for mixed space and time derivatives.
Combining all of the above, it follows that $\partial_t^m u(t) \in C\left( [0,T_*),\Sigma \right).$
\qed

\section{Derivation of invariants}\label{app:A}

The following invariants for \eqref{eq:inhomalb1} were introduced in \cite{Ribal2013AlberEquation}. We include them here, along with their derivation in our notation, for completeness.

\begin{lemma}\label{lm:inva} Let $u$  satisfy equation \eqref{eq:inhomalb1} with Hermitian initial data, and denote
\begin{align}
I_1[u] &:= \int\limits_x u(x,y,t)\Big|_{y=x} dx , \\
I_2[u] &:=  \int\limits_x (\partial_x - \partial_y)u(x,y,t)\Big|_{y=x} dx , \\
I_3[u] &:= \left( \frac{q}p \int\limits_x u^2(x,y,t)\Big|_{y=x} dx + \int\limits_x (\partial_x - \partial_y)^2u(x,y,t)\Big|_{y=x} dx \right) .
\end{align}
 Then, 
\begin{equation}\label{eq:lmstconcs12}
\frac{d}{dt} I_j[u] = 0 \qquad \mbox{ for } j\in \{1,2,3\}.	
\end{equation}
	
\end{lemma}

\medskip

\noindent {\bf Remark:} The function $u(x,x,t)$ is real-valued by the symmetry of the problem, but not non-negative. Thus the conservation of $I_1[u]$ is compatible with potentially large values of $u,$ but there will have to be positive and negative values so that the total integral is preserved. 

\medskip

\noindent {\bf Proof:} To prove this lemma, it is helpful to use the
change of variables
\begin{equation}
	\left\{
	\begin{array}{c}
		X=\frac{x+y}2 \\
		r = x-y
	\end{array}
	\right\} 
	\iff 
	\left\{
	\begin{array}{c}
	x = X+\frac{r}2 \\
	y = X - \frac{r}2
	\end{array}
	\right\}.
\end{equation}
Originally $x,y$ are two position variables,  and $X$ has the physical meaning of ``mean position'' while $r$ of ``lag'' between the two positions.  For any $j=1,2,\dots,$ observe that $\partial_x^j = (\partial_X + \frac{1}2\partial_r)^j,$ $\partial_y^j = (\partial_X - \frac{1}2\partial_r)^j.$  

Now denote
\begin{equation}\label{eq:chngvra}
	f=f(X,r,t) = u(x,y,t), \qquad n(X,t):= f(X,0,t).
\end{equation}
So if $u(x,y,t)$ satisfies equation \eqref{eq:inhomalb1},
then $f$ satisfies
\begin{equation}\label{eq:alberXr}
	i f_t  - p \partial_X \cdot \partial_r f - q \Bigl(n(X+\frac{r}2)-n(X-\frac{r}2)\Bigr) \, \Bigl( 
	\Gamma(r) + f 	\Bigr)=0.
\end{equation}

In this formalism, equation \eqref{eq:lmstconcs12} is elaborated to
\begin{align}
	&\frac{d}{dt} \int\limits_X n(X,t) dX = 0, \label{eq:I1consb} \\
	&\frac{d}{dt} \int\limits_X \partial_r f(X,r,t)\Big|_{r=0} dX = 0, \label{eq:I2consb}\\
	&\frac{d}{dt} \left( \frac{q}p \int\limits_X n^2(X,t) dX + \int\limits_X \partial_r^2 f(X,r,t)\Big|_{r=0} dX \right) = 0 \label{eq:I3consb}
\end{align}

\medskip 

Proof of equation \eqref{eq:I1consb}: setting $r=0$ in equation \eqref{eq:chngvra} and integrating in $X$ leads to
\[
\int\limits_X\Bigl( if_t(X,0,t) - p \partial_X \partial_r f(X,0,t) \Bigr)dX = 0 \implies \frac{d}{dt} \int\limits_X f(X,0,t) dX =0
\]
by performing integration by parts in $X.$

\medskip 

Proof of equation \eqref{eq:I2consb}: differentiating equation \eqref{eq:alberXr} with respect to $r_j$ yields
\begin{equation}\label{eq:derralb}
\begin{array}{c}
		i \partial_{t}f_{r_j}  - p \partial_X \cdot \partial_r f_{r_j} - \frac{q}2 \Bigl(n_{X_j}(X+\frac{r}2,t)+n_{X_j}(X-\frac{r}2,t)\Bigr) \, \Bigl( 
	\Gamma(r) + f 	\Bigr)  \hfill \\
	\hfill - q \Bigl(n(X+\frac{r}2)-n(X-\frac{r}2)\Bigr) \, \Bigl( 
	\partial_{r_j}\Gamma(r) + f_{r_j} 	\Bigr) =0.
\end{array}
\end{equation}
Now setting $r=0$ this becomes
\[
i\,\partial_t f_{r_j}(X,0,t)-p\,\partial_X\!\cdot\!\partial_r f_{r_j}(X,0,t)
-q\bigl(\Gamma(0)+n(X,t)\bigr)\,n_{X_j}(X,t)=0.
\]
Integrating over \(X\)  gives
\[
\begin{array}{c}
i\frac{\mathrm d}{\mathrm dt}\int\limits_X f_{r_j}(X,0,t)\,\mathrm dX
-q\int\limits_X\!\bigl(\Gamma(0)+n(X,t)\bigr)\,n_{X_j}(X,0,t)\,\mathrm dX=0, \\[8pt]

i\frac{\mathrm d}{\mathrm dt}\int\limits_X f_{r_j}(X,0,t)\,\mathrm dX
-q\Gamma(0)\int\limits_X\,n_{X_j}(X,t)\,\mathrm dX +\frac{q}2\int\limits_X\,\partial_{X_j}n^2(X,t)\,\mathrm dX=0. \\
\end{array}
\]
Equation \eqref{eq:I2consb} follows by  integration by parts in the two last terms.

\medskip 

Proof of equation \eqref{eq:I3consb}: by differentiating twice in $r,$ equation \eqref{eq:alberXr} becomes
\[
\begin{array}{r}
i\partial_t \Delta_r f - p \partial_X \cdot \partial_r \Delta_r f - \frac{q}4 \Bigl( \Delta_X n(X+\frac{r}2) - \Delta_X n(X-\frac{r}2) \Bigr) \Bigl( \Gamma(r) + f(X,r,t) \bigr)  \\
-q \nabla_X\Bigl( n(X+\frac{X}2) + n(X-\frac{n}2) \Bigr) \cdot \nabla_r \Bigl( \Gamma(r) + f(X,r,t) \Bigr) \\
-q \Bigl( n(X+\frac{r}{2}) - n(X-\frac{r}2) \Bigr) \Delta_r\Bigl( \Gamma(r) + f(X,r,t) \Bigr)=0.
\end{array}
\]
By setting $r=0$  we get
\[
\begin{array}{r}
i\partial_t \Delta_r f(X,0,t) - p \partial_X \cdot \partial_r \Delta_r f(X,0,t) \hfill \\
-q \nabla_X\Bigl( n(X+\frac{r}2) + n(X-\frac{r}2) \Bigr) \cdot \nabla_r \Bigl( \Gamma(r) + f(X,0,t) \Bigr) =0;
\end{array}
\]
by integrating in $X$ and some obvious integrations by parts
\begin{equation}\label{eq:last3a}
\begin{array}{r}
i\frac{d}{dt}\int\limits_X \Delta_r f(X,0,t) dX 
-2q \int\limits_X \nabla_X n(X,t)  \cdot  \nabla_rf(X,0,t)  dX =0.
\end{array}
\end{equation}
Moreover, by multiplying equation \eqref{eq:alberXr} with $f,$  setting $r=0$ and integrating in $X$ we get
\begin{equation}\label{eq:last23a}
\begin{array}{c}
\int\limits_X f(X,0,t) \Bigl( 
	i f_t(X,0,t)  - p \partial_X \cdot \partial_r f(X,0,t) \Bigr) dX=0 \implies \hfill \\
	\hfill \implies -\frac{i}{p} \frac{1}2 \frac{d}{dt} \int\limits_X n^2(X)dX =  \int\limits_X \nabla_X n(X,t) \cdot \nabla_r f(X,0,t) dX
\end{array}
\end{equation}
By combining equations \eqref{eq:last3a} and \eqref{eq:last23a}, relation \eqref{eq:I3consb} follows.
\qed

\end{document}